\documentclass{iopart}

\usepackage{iopams}

\expandafter\let\csname equation*\endcsname\relax

\expandafter\let\csname endequation*\endcsname\relax

\usepackage{amsmath,amssymb,amsthm}
\usepackage[english]{babel}
\usepackage[utf8x]{inputenc}
\usepackage[colorlinks=true,linkcolor=black,citecolor=blue]{hyperref}
\usepackage[capitalise,noabbrev]{cleveref}
\usepackage{pgfplots}
\pgfplotsset{compat=1.13}
\pgfplotsset{
  table/search path={Plotdata},
}

\usepackage[firstpage=true,color=black,placement=top,scale=1,opacity=1,vshift=-1cm]{background}

\SetBgContents{
\begin{minipage}{\linewidth}
\small
A definitive version of this manuscript was published in \textbf{Inverse Problems, Volume 35, Number 5, 2019} and is available at \url{https://dx.doi.org/10.1088/1361-6420/ab0b04}
\end{minipage}
}

\makeatletter
\renewcommand\@appendixstar{\@@par
 \ifnumbysec 
 \@addtoreset{table}{section}
 \@addtoreset{figure}{section}\fi
 \setcounter{section}{0}
 \setcounter{subsection}{0}
 \setcounter{subsubsection}{0}
 \setcounter{equation}{0}
 \setcounter{figure}{0}
 \setcounter{table}{0}
 \def\thesection{\Alph{section}} 
 \def\theequation{\ifnumbysec
      \Alph{section}.\arabic{equation}\else
      \Alph{section}\arabic{equation}\fi}
 \def\thetable{\ifnumbysec
      \Alph{section}\arabic{table}\else
      A\arabic{table}\fi}
 \def\thefigure{\ifnumbysec
      \Alph{section}\arabic{figure}\else
      A\arabic{figure}\fi}}
\makeatother

\newcommand{\N}{\mathbb{N}}
\newcommand{\NN}{\mathbb{N}_{0}}
\newcommand{\Z}{\mathbb{Z}}

\newcommand{\R}{\mathbb{R}}
\newcommand{\C}{\mathbb{C}}

\newcommand{\ii}{\text{i}}

\newcommand{\dx}{\,\mathrm{d}x}
\newcommand{\dy}{\,\mathrm{d}y}
\newcommand{\dz}{\,\mathrm{d}z}

\newcommand{\abs}[1]{\left| #1 \right|}
\newcommand{\norm}[1]{\left|\hspace*{-0.1em}\left| #1 \right|\hspace*{-0.1em}\right|}
\newcommand{\nnorm}[1]{\big|\hspace*{-0.1em}\big| #1 \big|\hspace*{-0.1em}\big|}

\renewcommand{\Re}{\mathrm{Re}}

\DeclareMathOperator{\supp}{supp}

\DeclareMathOperator{\sinc}{sinc}
\DeclareMathOperator{\Si}{Si}

\newtheorem{theorem}{Theorem}[section]
\newtheorem{proposition}[theorem]{Proposition}
\newtheorem{lemma}[theorem]{Lemma}
\newtheorem{corollary}[theorem]{Corollary}
\newtheorem{definition}[theorem]{Definition}

\graphicspath{{./Images/}}

\begin{document}

\title[Joint exit wave reconstruction and image registration]{A least-squares functional for joint exit wave reconstruction and image registration}
\author{Christian Doberstein and Benjamin Berkels}
\address{AICES Graduate School, RWTH Aachen University, Germany}

  \begin{abstract}
    Images generated by a transmission electron microscope (TEM)
    are blurred by aberrations from the objective lens and can
    be difficult to interpret correctly. One possible solution to this problem is to
    reconstruct the so-called exit wave, i.e.\ the electron wave in the microscope right
    before it passes the objective lens, from a series of TEM images acquired with varying
    focus. While the forward model of simulating a TEM image from a given exit wave is
    known and easy to evaluate, it is in general not possible to reconstruct the exit wave
    from a series of images analytically. The corresponding inverse problem can be
    formulated as a minimization problem, which is done in the well known MAL and MIMAP
    methods. We propose a generalization of these
    methods by performing the exit wave reconstruction and the registration of the image
    series simultaneously. We show that our objective functional is not convex with
    respect to the exit wave, which also carries over to the MAL and MIMAP
    functionals. The main result is the existence of minimizers of our objective
    functional. These results are based on the properties of a generalization of the
    cross-correlation. Finally, the applicability of our method is verified with a
    numerical experiment on simulated input data.
  \end{abstract}

  \section{Introduction}
    In transmission electron microscopy (TEM), an image is generated by recording a beam
    of electrons that was transmitted through a sample, also called specimen. Due to
    diffraction inside the specimen, the electrons carry structural information when
    leaving the specimen at its exit plane. The electrons are then diverted by
    electromagnetic lenses to form an image in the image plane of the electron microscope.
    
    There are several difficulties associated with directly interpreting TEM
    images. For instance, depending on the particular value of the objective lens
    aberrations at the time of image acquisition, entire atomic columns may not be visible
    in the recorded image or other features in the image may be mistakenly interpreted as
    atoms. These are the three main limitations inherent to TEM imaging:
    \begin{enumerate}
      \item The images are blurred by aberrations from the objective lens;
      \item temporal and spatial partial coherence of the illumination; and
      \item the phase of the image plane electron wave is missing in the images, since
            the microscope's camera only records the squared amplitude of the electron
            wave.
    \end{enumerate}
    Several approaches have been developed to deal with these limitations. These can be
    divided into two groups: improving the microscope's hardware components
    \cite{erni15, kabius09, haider98, otten93} and processing the TEM images
    \cite{kirkland16, beleggia04, kirkland04, coene96, opdebeeck96}.
    
    One of the approaches on the image processing side is exit wave reconstruction, which
    aims to reconstruct the electron wave at the exit plane of the specimen, the
    so-called exit wave. The exit wave is most commonly reconstructed from a series of
    TEM images taken with varying focus of the objective lens; this way each of the
    real-valued images encodes a different portion of the complex-valued exit wave.
    Although the reconstructed exit wave is theoretically free from aberrations of the
    objective lens, in practice it is frequently necessary to correct residual
    aberrations \cite{ophus16}. This is due to the fact that the aberration coefficients
    are oftentimes not known to a high precision.
    
    In this article, we investigate a
    generalization of two well-established variational methods for exit wave
    reconstruction: the multiple input maximum a-posteriori (MIMAP) algorithm
    \cite{kirkland84} and the maximum-likelihood (MAL) algorithm \cite{coene96}. In both
    algorithms, the exit wave is reconstructed by minimizing a functional of the kind
    \begin{equation*}
      J_\alpha[\Psi] = \frac{1}{N}\sum_{j=1}^N \norm{\mathcal{I}_{\Psi,Z_j}^\text{sim} - \mathcal{I}_j^\text{exp}}_{L^2}^2 + \alpha\norm{\Psi - \Psi_M}_{L^2}^2, \qquad{\alpha\ge0}.
    \end{equation*}
    It takes as input a finite series of experimental images
    $(\mathcal{I}_j^{\text{exp}})_{j=1,\ldots,N}$, which have been acquired with
    varying focus $(Z_j)_{j=1,\ldots,N}$ of the objective lens. The experimental images
    are compared to simulated images
    $(\mathcal{I}_{\Psi,Z_j}^{\text{sim}})_{j=1,\ldots,N}$ that are calculated
    from an estimate of the exit wave $\Psi$ and the corresponding focus values $Z_j$. The
    parameter $\Psi_M$ is an a-priori estimate of the exit wave and is only used in the
    MIMAP algorithm, whereas $\alpha$ is equal to zero in the MAL algorithm. Minimizing
    the functional $J_\alpha$ corresponds to finding a (regularized) least squares
    approximation to the experimental input images. See, for example, \cite{barthel10} for
    an application of the MAL algorithm.
    
    In principle, the exit wave can be reconstructed by only minimizing $J_\alpha$. As an
    initial guess, a constant exit wave corresponding to the square root of the mean image
    intensities or the results from one iteration of the paraboloid method may be used
    \cite{coene96,opdebeeck96}. However, due to specimen drift and microscope
    instabilities, an accurate registration of the experimental image series is necessary
    in addition to the minimization of $J_\alpha$. In the MAL algorithm as described in
    \cite{coene96}, this problem is solved in two steps. First, a rough initial alignment
    of the experimental images is computed by calculating the cross-correlation of
    subsequent images in the series. Second, after every iteration of the minimization
    method the alignment is further improved by calculating the cross-correlation of each
    experimental image with the corresponding simulated image. In the MIMAP
    algorithm the alignment is improved in every iteration by additionally applying a
    translation to the simulated images $\mathcal{I}_{\Psi, Z_j}^\text{sim}$ in the
    functional and minimizing the functional's value with respect to this translation
    \cite{kirkland84}. This is done without taking the update of the exit wave in the
    current iteration into account.
    
    The generalization of the MAL and MIMAP algorithms we propose is to
    treat this as joint optimization problem, i.e.\ to optimize the exit wave and the registration
    simultaneously rather than alternatingly. The generalized functional is
    \begin{equation*}
      E_\alpha[\Psi,t_1,\ldots,t_N] = \frac{1}{N}\sum_{j=1}^N \norm{\mathcal{I}_{\Psi,Z_j}^{\text{sim}} - \mathcal{I}_{t_j,j}^{\text{exp}}}_{L^2}^2 + \alpha\norm{\Psi - \Psi_M}_{L^2}^2 \qquad{\alpha\ge0},
    \end{equation*}
    where $\mathcal{I}_{t_j,j}^{\text{exp}}$ is the experimental image
    $\mathcal{I}_j^{\text{exp}}$ translated by $t_j\in\R^2$.
    This coupled treatment of the reconstruction with the registration
    is essential to our analysis of the existence of minimizers.
    
    Carrying out such a mathematical analysis is the main contribution of this paper.
    To this end, we show that our objective
    functional is weakly lower semi-continuous and coercive, which enables us to proof
    the existence of minimizers with the direct method. Additionally, we show that our
    objective functional is in general not convex with respect to the exit wave. These
    results are fundamentally based on 1) a novel factorization property of the
    weight used in the forward model, which is given in \cref{tcc_factorization}, and 2) the notion of the weighted cross-correlation
    developed in \cref{wcc_appendix}.
    
    \medskip\noindent The structure of this article is as follows:
    \begin{itemize}
      \item \textit{The forward model:}\quad The simulation of a TEM image from a given
            exit wave amounts to calculating the weighted autocorrelation of the exit
            wave. In \cref{forward_model}, the weight is introduced in detail and some of
            its properties are presented. The definition of the weighted
            cross-correlation is given in \cref{wcc_appendix}, together with several
            generic properties independent of the application to exit wave
            reconstruction.
      \item \textit{The inverse problem:}\quad In \cref{objective_functional}, our
            objective functional for exit wave reconstruction and its derivatives are
            given. It is also shown that the functional is not convex with respect to the
            exit wave unless all experimental input images are zero.
      \item \textit{Existence of minimizers:}\quad The results from
            \cref{lowpass_filter_appendix} suggest that the objective functional as given
            in \cref{objective_functional} is not coercive. This problem is solved by
            adding a Tikhonov regularization term to the functional. Using the direct
            method, it is shown in \cref{existence} that minimizers of the regularized
            functional exist.
      \item \textit{Numerical experiment:}\quad We conclude with a numerical experiment
            on synthetic input data in \cref{experiment}, showing that our objective
            functional can indeed be used to reconstruct the exit wave and register the
            images simultaneously in practice.
    \end{itemize}
    
    \medskip\noindent Although this article is
    written from the point of view of exit wave reconstruction, the results are more
    generally applicable to any least-squares inverse problem where the forward model is a
    weighted autocorrelation as in \cref{forward_model}. More precisely, the existence and
    convexity results still hold if the specific weight used for TEM image simulation is
    replaced with any weight that satisfies
    \cref{elementary_tcc_properties,tcc_continuity,tcc_factorization}.
    
    In the following, we frequently use the Fourier transform. Here, we use
    \begin{equation*}
      \big(\mathcal{F}f\big)(x) := \int_{\R^d} f(y)e^{-2\pi\ii x\cdot y}\dy \tag*{$\forall\,x\in\R^d$}
    \end{equation*}
    as the definition of the Fourier transform of an integrable function
    $f\in L^1(\R^d, \C)$. This definition is convenient, since no additional scaling
    factor is needed when applying the convolution theorem. If the integration domain is
    omitted from an integral expression, the domain is always the entire vector space
    $\R^d$. Furthermore, all considered subsets of $\R^d$ are implicitly assumed to be
    Lebesgue measurable.
   
  \section{The forward model: simulating TEM images}\label{forward_model}
    The simulation of TEM images from a given exit wave and focus value is one of the
    core components of the inverse problem. In Fourier space, the simulated image
    $\mathcal{I}_{\Psi,Z}^{\text{sim}}\in L^2(\R^2, \C)$ is given by a weighted
    autocorrelation of the exit wave $\Psi\in L^2(\R^2, \C)$ \cite{kirkland10}.
    Explicitly, we have
    \begin{equation}\label{image_simulation_eq}
      \mathcal{I}_{\Psi,Z}^{\text{sim}} = \Psi \star_{\widehat{T}_Z} \Psi,
    \end{equation}
    where $Z\in\R$ is the focus value, $\widehat{T}_Z:\R^2 \times \R^2 \rightarrow \C$ is
    the transmission cross-coefficient (TCC), and $\Psi \star_{\widehat{T}_Z} \Psi$
    denotes the weighted cross-correlation (cf.\ \cref{wcc_def}). Taking partial
    coherence into account, the TCC is defined as
    \begin{equation}\label{general_tcc}
      \widehat{T}_Z: \R^2 \times \R^2 \rightarrow \C, \; (v,w) \mapsto \int_{\R}\int_{\R^2} s(u) f(Z') t_{Z+Z'}(v+u)t_{Z+Z'}^*(w+u)\,\mathrm{d}u\,\mathrm{d}Z',
    \end{equation}
    where $s\in L^1(\R^2,\R)$ and $f\in L^1(\R,\R)$ are probability density functions and
    $t_Z:\R^2\to\C$ is the pupil function. Here, $s$ is the normalized intensity
    distribution of the illumination and $f$ is the normalized focus spread. The pupil
    function is
    \begin{equation}
      t_Z(v) := p_Z(v)a(v) \tag*{$\forall\,v\in\R^2$,}
    \end{equation}
    where $p_Z:\R^2\to\C$ is the pure phase transfer function and $a:\R^2\to\R$ is the
    aperture function.
    
    The pure phase transfer function models the aberrations of the objective lens and is
    given by
    \begin{equation}
      p_Z(v) := \exp\big(-2\pi\ii\chi_Z(v)\big) \tag*{$\forall\,v\in\R^2$.}
    \end{equation}
    For the sake of simplicity, we only consider the two most important
    isotropic aberrations -- focus and third order spherical aberration -- in the wave
    aberration function $\chi_Z:\R^2\to\R$. Then
    \begin{equation}
      \chi_Z(v) := \frac{1}{2}Z\lambda\abs{v}^2 + \frac{1}{4}C_s\lambda^3\abs{v}^4 \tag*{$\forall\,v\in\R^2$,}
    \end{equation}
    where $\lambda\in\R_{>0}$ is the electron wavelength and $C_s\in\R$ is the
    coefficient of the spherical aberration. Additional terms for higher order or
    anisotropic aberrations can easily be included in the wave aberration function. The
    general formula including all possible aberrations can e.g.\ be found in
    \cite{kirkland10} or \cite{thust07}.
    
    The aperture function models the objective aperture and is given by
    \begin{equation}
      a(v) := \begin{cases}
               1, & \text{if $\lambda\norm{v}_2 < \alpha_\text{max}$,} \\
               0, & \text{otherwise,}
             \end{cases} \tag*{$\forall\,v\in\R^2$,}
    \end{equation}
    where $\lambda\in\R_{>0}$ is the electron wavelength and
    $\alpha_\text{max}\in\R_{>0}$ is the maximum semiangle allowed by the objective
    aperture. The set-theoretic support of the aperture function is a ball of radius
    $r_a := \alpha_\text{max}/\lambda$ centered at the origin and is denoted by
    $A := B_{\alpha_\text{max}/\lambda}(0)$.
    
    Note that the pure phase transfer function $p_Z$ is continuous and the aperture
    function $a$ is bounded. The parameters $\lambda$,
    $C_s$ and $\alpha_\text{max}$ are treated as constants in the following.
    
    \medskip
    Depending on the particular probability densities $s$ and $f$, the support of the TCC
    $\widehat{T}_Z$ might be unbounded. As a TCC with bounded support significantly
    simplifies the theory developed here, we will consider
    \begin{equation}\label{special_tcc}
      T_Z(v,w) := a(v)a^*(w) \int_\R \int_{\R^2} s(u)f(Z')p_{Z+Z'}(v+u)p_{Z+Z'}^*(w+u)\,\mathrm{d}u\,\mathrm{d}Z'
    \end{equation}
    instead of $\widehat{T}_Z$ in the following. In doing so, we ignored the dependency
    of the aperture function on the integration variable $u$. This is a common
    simplification that is justified if the diameter of the aperture is sufficiently
    large compared to the highest frequency of interest
    \cite{kirkland10,ishizuka80,frank73}. In particular, this can be considered as the
    first step towards Ishizuka's well established approximation to the TCC
    $\widehat{T}_Z$ \cite{ishizuka80}, which is used in the MAL and MIMAP
    algorithms.
    
    \begin{lemma}[Elementary properties of the TCC]\label{elementary_tcc_properties}
      The transmission cross-coefficient $T_Z$ has the following properties:
      \begin{enumerate}
        \item $\abs{T_Z(v,w)}\le 1$ for all $v,w\in\R^2$,
        \item $T_Z^*(v,w) = T_Z(w,v)$ for all $v,w\in\R^2$,
        \item $T_Z(v,v) = 1$ for all $v\in A$,
        \item $T_Z(v,w) = 0$ for all $v,w\in\R^2\backslash A$.
      \end{enumerate}
    \end{lemma}
    \noindent The second property has also
    been observed by Ishizuka in \cite{ishizuka80}.

    \proof All of the properties follow immediately from the definition of $T_Z$ and
      the properties of its components.\qed
    
    \medskip
    The first two properties of the TCC in \cref{elementary_tcc_properties} also hold
    for the more general version $\widehat{T}_Z$ of the TCC. In \cite{ishizuka80}, it is
    shown that, as a direct consequence of the second property, the Fourier space image
    $G := \mathcal{I}_{\Psi,Z}^{\text{sim}}$ satisfies Friedel's law, i.e.
    \begin{equation}
      G(v) = G^*(-v) \tag*{$\forall\,v\in\R^2$.}
    \end{equation}
    Taking the inverse Fourier transform, this shows that the simulated image in real
    space is indeed real-valued, as one would expect from a TEM image.
    
    \medskip
    An important property of the ordinary cross-correlation is that $f\star g$ is
    continuous for functions $f\in L^p(\R^d)$ and $g\in L^q(\R^d)$, where $d\in\N$ and
    $p,q>1$ with $\frac{1}{p} + \frac{1}{q} = 1$. This property is generalized to
    weighted cross-correlations in \cref{wcc_continuity} on the additional assumption
    that the weight itself is continuous on a suitable open subset of $\R^d$.
    
    \begin{lemma}[Continuity of the TCC]\label{tcc_continuity}
      The transmission cross-coefficient $T_Z$ is continuous on $A \times A$.
    \end{lemma}
    \proof The integrand of $T_Z$ is dominated by the function
      $g\in L^1(\R^2\times\R, \R)$ defined as $g(u,Z') := \abs{s(u)f(Z')}$. Therefore,
      the result follows from the dominated convergence theorem and the continuity of
      $p_Z$.\qed
    
    \medskip
    If $(s_n)_{n\in\N}\in L^1(\R^2,\R)^\N$ and $(f_n)_{n\in\N}\in L^1(\R,\R)^\N$ are
    sequences of probability density functions with uniformly bounded support such that
    \begin{equation*}
      \lim_{n\rightarrow\infty} \int_{B_\varepsilon(0)} s_n(u)\,\mathrm{d}u = 1
      \qquad\text{and}\qquad
      \lim_{n\rightarrow\infty} \int_{(-\varepsilon,\varepsilon)} f_n(Z')\,\mathrm{d}Z' = 1
    \end{equation*}
    for all $\varepsilon>0$, then
    \begin{align*}
      \lim_{n\rightarrow\infty}\lim_{m\rightarrow\infty}&a(v)a^*(w)\int_{\R}\int_{\R^2} s_n(u) f_m(Z') p_{Z+Z'}(v+u)p_{Z+Z'}^*(w+u)\,\mathrm{d}u\,\mathrm{d}Z' \\
      =\;&a(v)a^*(w)p_Z(v)p_Z^*(w) = t_Z(v)t_Z^*(w)
    \end{align*}
    for all $v,w\in\R^2$. This corresponds to the limiting case of perfect spatial and
    temporal coherence, in which the TCC is simply a product of the pupil function and
    its complex conjugate. In particular, the image simulation simplifies to
    \begin{equation*}
      \mathcal{I}_{\Psi,Z}^{\text{sim}} = \Psi \star_{T_Z} \Psi
      = (\Psi t_Z) \star (\Psi t_Z)
      = (\Psi p_Z a) \star (\Psi p_Z a).
    \end{equation*}
    
    In the general case, there is no straightforward way to factorize the TCC and
    express $\Psi \star_{T_Z} \Psi$ as an ordinary cross-correlation. However, it is
    possible to approximate $T_Z$ by finite sums of factorizable functions, which is made
    precise in the following proposition. This factorization property will be most
    helpful for the analysis of the functional, as it enables us to apply the convolution
    theorem to weighted cross-correlations.
    
    In the following, the integrand of the TCC is split into the two parts
    \begin{align*}
      g&: \R^2 \times \R \rightarrow \R, \quad (u,Z') \mapsto s(u)f(Z'), \\
      q_{v,w}&: \R^2 \times \R \rightarrow \C, \quad (u,Z') \mapsto p_{Z+Z'}(v+u)p_{Z+Z'}^*(w+u)
    \end{align*}
    so that $T_Z(v,w)=a(v)a^*(w)\iint g(u,Z')q_{v,w}(u,Z')\,\mathrm{d}u\,\mathrm{d}Z'$
    for all $v,w \in \R^2$. Here, the supremum norm of a vector-valued function
    $\alpha: \R^n\rightarrow\C^m$ is defined as
    \begin{equation*}
      \norm{\alpha}_\infty := \sup_{x\in\R^n} \norm{\alpha(x)}_\infty.
    \end{equation*}
    
    \begin{proposition}[Factorization property of the TCC]\label{tcc_factorization}
      Assume that $s$ and $f$ are differentiable with
      \begin{equation*}
       \sup_{v,w\in A}\norm{g\nabla q_{v,w}}_\infty<\infty\qquad\text{and}\qquad\norm{\nabla g}_\infty<\infty.
      \end{equation*}
      Then there exists a sequence $(T_{Z,N})_{N\in\N}$ of functions converging uniformly
      to $T_Z$ with the following properties:
      \begin{enumerate}
        \item $T_{Z,N}$ is continuous and bounded on $A\times A$ and zero on
              $(\R^2\times\R^2)\backslash (A\times A)$ for all $N\in\N$.
        \item For all $N\in\N$, $j\in\{1,\ldots,N\}$ there are functions
              $t_{Z,N,j}:\R^2\rightarrow\C$ such that $t_{Z,N,j}$ is continuous and
              bounded on $A$, zero on $\R^2\backslash A$ and
              \begin{equation}
                T_{Z,N}(v,w) = \sum_{j=1}^N t_{Z,N,j}(v)t_{Z,N,j}^*(w) \tag*{$\forall\,v,w\in\R^2$.}
              \end{equation}
      \end{enumerate}
    \end{proposition}
    
    \proof Denote by
      \begin{equation*}
        T_Z'(v,w) := \int_\R \int_{\R^2} s(u)f(Z')p_{Z+Z'}(v+u)p_{Z+Z'}^*(w+u)\,\mathrm{d}u\,\mathrm{d}Z'
      \end{equation*}
      the TCC without the aperture function. As the integrand is continuous in $(u,Z')$ and already factorized with
      respect to $v$ and $w$, we expect the Riemann sums
      \begin{align*}
        &\sum_{\alpha=-M}^M \sum_{\beta,\gamma=-M}^M \underbrace{\delta_M \delta_M^2 s(u_{\beta,\gamma}) f(Z'_{\alpha})}_{\ge 0} p_{Z+Z'_{\alpha}}(v+u_{\beta,\gamma}) p_{Z+Z'_{\alpha}}^*(w+u_{\beta,\gamma}), \\
        &\bullet\; M\in\N, \\
        &\bullet\; \delta_M = M^{-k} \quad\text{for $k\in\big(\tfrac{3}{4},1\big)$}, \\
        &\bullet\; u_{\beta,\gamma} = (\beta\delta_M, \gamma\delta_M) \in \R^2, \\
        &\bullet\; Z_{\alpha}' = \alpha\delta_M \in\R
      \end{align*}
      to be suitable approximations of $T_Z'$ as sums of factorizable functions. If we define 
      \begin{equation}
        t_{Z,M,\alpha,\beta,\gamma}(v) := a(v)\sqrt{\delta_M \delta_M^2 s(u_{\beta,\gamma}) f(Z'_{\alpha})} p_{Z+Z'_{\alpha}}(v+u_{\beta,\gamma}) \tag*{$\forall\,v\in\R^2$,}
      \end{equation}
      then reordering the indices $\alpha,\beta,\gamma\in\{-M,\ldots,M\}$ to a single
      index $j\in\{1,\ldots,N\}$ yields the sought functions $t_{Z,N,j}$. It remains to
      show that the sums converge uniformly to $T_Z'$ on $A\times A$ as
      $M\rightarrow\infty$, since the aperture function is zero on $\R^2\backslash A$.
      
      \medskip
      Let $\varepsilon>0$ and fix $v,w\in A$. First, we note that the integrand
      \begin{equation*}
        h_{v,w}: \R^2 \times \R \rightarrow \C, \quad (u, Z') \mapsto s(u)f(Z')p_{Z+Z'}(v+u)p_{Z+Z'}^*(w+u)
      \end{equation*}
      is Lipschitz continuous. This follows from the fact that the gradient is bounded by
      \begin{align*}
        \norm{\nabla h_{v,w}}_\infty &= \norm{\nabla(g q_{v,w})}_\infty = \norm{g\nabla q_{v,w} + q_{v,w}\nabla g}_\infty
        \le \norm{g\nabla q_{v,w}}_\infty + \norm{\nabla g}_\infty \\
        &\le \sup_{v,w\in A}\norm{g\nabla q_{v,w}}_\infty + \norm{\nabla g}_\infty =: L < \infty.
      \end{align*}
      Let $J_M := [-M\delta_M,(M+1)\delta_M]$ and $I_x := [x\delta_M, (x+1)\delta_M]$ for
      all $x\in\Z$, $M\in\N$. Using the triangle inequality,
      \begin{align}
        &\phantom{=} \abs{\sum_{\alpha=-M}^M \sum_{\beta,\gamma=-M}^M \delta_M^3 h_{v,w}(u_{\beta,\gamma}, Z'_\alpha) - \int_{\R} \int_{\R^2} h_{v,w}(u,Z')\,\mathrm{d}u\,\mathrm{d}Z'} \notag\\
        &\le \abs{\sum_{\alpha=-M}^M \sum_{\beta,\gamma=-M}^M \delta_M^3 h_{v,w}(u_{\beta,\gamma}, Z'_\alpha) - \int_{J_M} \int_{J_M\times J_M} h_{v,w}(u,Z')\,\mathrm{d}u\,\mathrm{d}Z'} + C_M \notag\\
        &\le \sum_{\alpha=-M}^M \sum_{\beta,\gamma=-M}^M \int_{I_\alpha}\int_{I_\beta\times I_\gamma} \abs{h_{v,w}(u_{\beta,\gamma}, Z'_\alpha) - h_{v,w}(u,Z')}\,\mathrm{d}u\,\mathrm{d}Z' + C_M \notag\\
        &\le \sum_{\alpha=-M}^M \sum_{\beta,\gamma=-M}^M \int_{I_\alpha}\int_{I_\beta\times I_\gamma} 3L\norm{(u_{\beta,\gamma},Z_\alpha') - (u,Z')}_\infty\,\mathrm{d}u\,\mathrm{d}Z' + C_M \notag\\
        &\le (2M+1)^3\delta_M^3 \cdot 3L\delta_M + C_M, \label{tcc_factorization_upperbound}
      \end{align}
      where $C_M = \sup_{v,w\in A} \int_{\R^3\backslash J_M^3} \abs{h_{v,w}(x)}\dx$. The
      functions
      \begin{equation*}
        c_M:\overline{A}\times\overline{A} \rightarrow \R, \quad (v,w)\mapsto\int_{\R^3\backslash J_M^3} \abs{h_{v,w}(x)}\dx \qquad (M\in\N)
      \end{equation*}
      define a monotone sequence of continuous functions, which converges pointwise to
      zero on the compact space $\overline{A}\times\overline{A}$. Hence $c_M$ converges
      uniformly to zero as $M\rightarrow\infty$ by Dini's theorem, which implies
      $\lim_{M\rightarrow\infty} C_M = 0$. The first summand in
      \cref{tcc_factorization_upperbound} also converges to zero for
      $M\rightarrow\infty$, since $\delta_M = M^{-k}$ with $k\in(\frac{3}{4},1)$. This
      shows the uniform convergence of the Riemann sums to $T_Z'$ on $A \times A$.\qed
    
    \medskip
    Intuitively, the assumption $\sup_{v,w\in A}\norm{g\nabla q_{v,w}}_\infty<\infty$
    ensures that $s$ and $f$ decrease at least as fast as the gradient of $q_{v,w}$
    grows. For instance, if $s$ and $f$ are Gaussian probability densities, then this
    condition is clearly fulfilled as $\nabla q_{v,w}$ grows only polynomially.
    
    \medskip
    Using the factorization property of $T_Z$, it is shown in the following corollary
    that $\mathcal{F}^{-1}(\mathcal{I}_{\Psi,Z}^{\text{sim}})$ is not only real-valued,
    but also nonnegative. This is clearly a desired property of the simulation, since TEM
    images consist of electron counts and thus should not contain negative values.
    
    \begin{corollary}\label{tcc_nonnegative}
      The simulated real space images $\mathcal{F}^{-1}(\mathcal{I}_{\Psi,Z}^{\text{sim}}) = \mathcal{F}^{-1}(\Psi\star_{T_Z}\Psi)$
      are nonnegative for all $\Psi\in L^2(\R^2, \C)$.
    \end{corollary}
    \proof By \cref{tcc_factorization}, the TCC is an element of $W^+(\R^2)$
      (cf.\ \cref{wcc_appendix}). Combining the elementary properties from
      \cref{elementary_tcc_properties} with \cref{auto_wcc_nonnegative} therefore shows
      that $\mathcal{F}^{-1}(\Psi\star_{T_Z}\Psi)\ge 0$ for all $\Psi\in L^2(\R^2, \C)$.
      \qed
    
    \medskip
    For the remainder of this section, we
    describe how the forward model from the MAL and MIMAP algorithms fits into the
    framework described here. There, an image is simulated using
    \cref{image_simulation_eq} with Ishizuka's TCC \cite{ishizuka80} instead of
    $\widehat{T}_{Z}$ or $T_{Z}$. Ishizuka's TCC is an approximation to $T_{Z}$ that makes
    two assumptions for numerical efficiency:
    \begin{enumerate}
      \item the use of a high resolution TEM to investigate a weakly scattering object;
      \item the intensity distribution of the illumination $s$ and the focus spread $f$
            can be modeled by Gaussian probability densities.
    \end{enumerate}
    Using these assumptions, the TCC $T_{Z}$ is simplified to yield
    \begin{equation}\label{MAL_tcc}
      T_{Z}^\text{Ishizuka}(v,w) := p_Z(v)p_Z^*(w)a(v)a^*(w)E_s(v,w)E_t(v,w),
    \end{equation}
    where $E_s$ and $E_t$ are damping envelopes that model the effects of partial spatial
    and temporal coherence respectively. They are given by
    \begin{align*}
      E_s(v,w) &:= \exp\left( -\left(\frac{\pi\alpha}{\lambda}\right)^2 \abs{\nabla\chi(v) - \nabla\chi(w)}^2 \right), \\
      E_t(v,w) &:= \exp\left( -\frac{1}{2}(\pi\Delta\lambda)^2 \left(\abs{v}^2 - \abs{w}^2\right)^2 \right),
    \end{align*}
    where $\lambda\in\R_{>0}$ is the electron wavelength, $\alpha\in\R_{\ge 0}$ is the
    half angle of beam convergence and $\Delta\in\R_{\ge 0}$ is the focus spread
    parameter.
    
    The factorization property given in \cref{tcc_factorization} is related to the focal
    integration approximation that is used in the MAL algorithm. There, the temporal
    coherence envelope $E_t$ is approximated by a sum of factorizable terms that
    originate from a numerical integration of
    \begin{align}\label{tcc_mal_et_int}
      E_t(v,w) &= \int_{-\infty}^{\infty} f_{\Delta}(x) \exp\left(-\pi\ii\lambda x\left(\abs{v}^2-\abs{w}^2\right)\right)\dx, \\
      f_{\Delta}(x) &= \frac{1}{\sqrt{2\pi}\Delta}\exp\left(-\frac{x^2}{2\Delta^2}\right), \notag
    \end{align}
    where $f_{\Delta}$ is the focus spread. This integral expression is equal to the
    original definition of $E_t$, since
    \begin{equation*}
      \exp(-C^2) = \int_{-\infty}^\infty \frac{1}{\sqrt{\pi}}\exp(-2\ii Cx)\exp(-x^2)\dx
    \end{equation*}
    holds for all $C\in\R$. For practical reasons, the spatial coherence
    envelope $E_s$ is only roughly approximated by $E_s(v,w) \approx E_s(v,0) E_s^*(w,0)$
    in the MAL algorithm. Therefore, the approximations to the TCC used in the MAL
    algorithm do not converge to $T_Z^\text{Ishizuka}$. However, if the
    spatial coherence envelope is instead expressed as
    \begin{align}\label{tcc_mal_es_int}
      E_{s}(v,w) &= \prod_{j=1}^2 \exp\left( -\left(\frac{\pi\alpha}{\lambda}\right)^2 \abs{\partial_j\chi(v) - \partial_j\chi(w)}^2 \right) \notag\\
      &= \prod_{j=1}^2 \int_{-\infty}^\infty \frac{1}{\sqrt{\pi}\alpha}\exp\left(-\frac{x^2}{\alpha^2}\right)
      \exp\left(-\frac{2\pi ix}{\lambda}\left(\partial_j\chi_Z(v) - \partial_j\chi_Z(w)\right)\right)\dx,
    \end{align}
    then a numerical integration of \cref{tcc_mal_et_int,tcc_mal_es_int} yields
    approximations to $T_Z^\text{Ishizuka}$ as sums of factorizable functions. It
    can be shown that these approximations converge uniformly, so the approximation
    property in \cref{tcc_factorization} also holds for Ishizuka's TCC.
    
    From the definition of $T_Z^\text{Ishizuka}$ in \cref{MAL_tcc} it is clear that
    \cref{tcc_continuity} continues to hold as well as all of the elementary properties
    given in \cref{elementary_tcc_properties}.
  
  \section{The objective functional}\label{objective_functional}
  
  In this section, our objective functional for exit wave reconstruction and its
  derivatives are given, alongside with a result regarding the convexity of the
  objective functional. By the results of the previous section, the statements in this
  section still hold if the TCC $T_Z$ is replaced with $T_Z^\text{Ishizuka}$.
  
  \begin{definition}\label{functional_definition}
    Let $N\in\N$ and $g_1,\ldots,g_N\in L^1(\R^2,\R_{\ge 0}) \cap L^2(\R^2,\R_{\ge 0})$
    be a series of real space TEM images such that
    $\supp(\mathcal{F}(g_j))\subseteq \overline{2A}$ for all $j\in\{1,\ldots,N\}$. Denote by
    $Z_j\in\R$ the focus value associated with the image $g_j$. The objective functional
    for joint reconstruction and registration is defined as
    \begin{equation}\label{functional_definition_eq}
      E: L^2(A,\C) \times (\R^2)^N \rightarrow \R, \quad (\Psi, t) \mapsto \frac{1}{N}\sum_{j=1}^N \nnorm{\Psi\star_{T_{Z_j}}\Psi - \mathcal{F}(g_j\circ\phi_{t_j})}_{L^2}^2,
    \end{equation}
    where 
    $\phi_{y}: \R^2 \rightarrow \R^2, \; x \mapsto x + y$ is the translation by
    $y\in\R^2$. Here,
    $I_{\Psi,Z_j}^\text{sim} = \Psi\star_{T_{Z_j}}\Psi$ are the simulated images in
    Fourier space and $I_{t_j,j}^\text{exp} = \mathcal{F}(g_j\circ\phi_{t_j})$ are the
    translated experimental images in Fourier space.
  \end{definition}
  
  The regularization term $\alpha\norm{\Psi - \Psi_M}_{L^2}^2$ presented in the
  introduction is not relevant in this section and therefore will only be added to our
  objective functional in \cref{existence}.
  
  The functional is well-defined, because $\Psi\star_{T_{Z_j}}\Psi\in L^2(\R^2, \C)$ by
  \cref{wcc_lplp} and clearly also $\mathcal{F}(g_j\circ\phi_{t_j})\in L^2(\R^2, \C)$ for
  all $j\in\{1,\ldots,N\}$. The domain of the exit wave is restricted to $A$, since all
  frequencies outside of $A$ are filtered out due to
  \begin{equation}\label{objfunc_domain_restriction_eq}
    \Psi\star_{T_{Z_j}}\Psi = (\Psi a)\star_{T_{Z_j}}(\Psi a)
  \end{equation}
  and therefore do not contribute to the simulated images.
  \cref{objfunc_domain_restriction_eq} follows from the fact that the aperture function
  satisfies $a(v)^2 = a(v)$ for all $v\in\R^2$.
  The additional restriction
  $\supp(\mathcal{F}(g_j))\subseteq\overline{2A}$ on the frequencies of the experimental
  images does not affect the reconstruction, since
  $\supp(\Psi\star_{T_Z}\Psi)\subseteq\overline{2A}$ for all $\Psi\in L^2(A,\C), Z\in\R$
  and $\supp(\mathcal{F}(g_j\circ\phi_y)) = \supp(\mathcal{F}(g_j))$ for all $y\in\R^2$.
  Thus only image plane frequencies $v$ with $\abs{v} < 2r_a$, where $r_a$ is the aperture
  radius, can be reconstructed with the objective functional. Hence a low-pass filter of
  radius $2r_a$ may be applied to the experimental input images $g_1,\ldots,g_N$ without
  loss of generality.
  
  \medskip
  Note that a translation in real space can be expressed as a modulation in
  Fourier space by means of the identity
  \begin{equation*}
    \mu_y\mathcal{F}(g) = \mathcal{F}(g \circ \phi_y),
  \end{equation*}
  which holds for all $g\in L^1(\R^2, \R)$ and $y\in\R^2$, where
  $\mu_y: \R^2 \rightarrow \C, \; x \mapsto e^{2\pi i x\cdot y}$ is the modulation by $y$.
  Therefore, the objective functional in \cref{functional_definition_eq} can be written as
  \begin{equation*}
    E: L^2(A,\C) \times (\R^2)^N \rightarrow \R, \quad (\Psi, t) \mapsto \frac{1}{N}\sum_{j=1}^N \nnorm{\Psi\star_{T_{Z_j}}\Psi - \mu_{t_j}\mathcal{F}(g_j)}_{L^2}^2.
  \end{equation*}
  With this variant of the objective
  functional it can be seen that the translation may equivalently be applied to the exit
  wave instead of the experimental images. This follows from
  \begin{align*}
    \Big((\Psi\mu_y)\star_{T_{Z_j}}(\Psi\mu_y)\Big)(x) &= \int_{\R^2} \Psi^*(z)\mu_y^*(z)\Psi(x+z)\mu_y(x+z)T_{Z_j}(x+z,z)\dz \\
    &= \int_{\R^2}e^{2\pi i x\cdot y}\Psi^*(z)\Psi(x+z)T_{Z_j}(x+z,z)\dz \\
    &= \mu_y(x) \Big(\Psi\star_{T_{Z_j}}\Psi\Big)(x),
  \end{align*}
  which implies
  \begin{equation*}
    \frac{1}{N}\sum_{j=1}^N \nnorm{\Psi\star_{T_{Z_j}}\Psi - \mu_{t_j}\mathcal{F}(g_j)}_{L^2}^2 = \frac{1}{N}\sum_{j=1}^N \nnorm{(\Psi\mu_{-t_j})\star_{T_{Z_j}}(\Psi\mu_{-t_j}) - \mathcal{F}(g_j)}_{L^2}^2.
  \end{equation*}
  The right-hand side of the above
  equation is essentially identical to the functional that is used in the MIMAP algorithm
  for the optimization of the registration \cite{kirkland84}. Therefore, the MIMAP
  algorithm can be considered as a special case of our joint optimization approach.
  
  \medskip
  \cref{functional_definition_eq} can be regarded as the Fourier space formulation of
  $E$. Since the Fourier
  transform is unitary, the equivalent real space
  formulation of the objective functional is
  \begin{align*}
    E[\Psi,t] = \frac{1}{N}\sum_{j=1}^N \nnorm{\mathcal{F}^{-1}\big(\Psi\star_{T_{Z_j}}\Psi\big) - g_j \circ \phi_{t_j}}_{L^2}^2.
  \end{align*}
  
  In the MAL algorithm, the exit wave is reconstructed by alternatingly minimizing the
  MAL functional
  \begin{equation*}
    E_{\text{MAL}}[\Psi] = \frac{1}{N}\sum_{j=1}^N \nnorm{\Psi\star_{T_{Z_j}^\text{Ishizuka}}\Psi - \mathcal{F}(g_j)}_{L^2}^2
  \end{equation*}
  and updating the registration, where the registration is updated using the
  cross-correlation of simulated and experimental images. Interestingly,
  the MAL algorithm is equivalent to alternatingly minimizing
  our objective functional with respect to $\Psi$ and $t$, which
  can be seen as follows.
  
  On the one hand, if $t=0$ is fixed, then $E[\Psi,t] = E_\text{MAL}[\Psi]$ for all
  $\Psi\in L^2(A,\C)$ and it is clear that minimizing $E$ with respect to $\Psi$ is
  equivalent to minimizing $E_\text{MAL}$. On the other hand, if $\Psi$ is fixed, consider
  the real space formulation of $E$ and define
  $f_{\Psi,j} := \mathcal{F}^{-1}\big(\Psi\star_{T_{Z_j}}\Psi\big)$ for all
  $j\in\{1,\ldots,N\}$. Minimizing
  \begin{equation*}
    E[\Psi,t] = \frac{1}{N}\sum_{j=1}^N \norm{f_{\Psi,j} - g_j \circ \phi_{t_j}}_{L^2}^2
    = \frac{1}{N}\sum_{j=1}^N \norm{f_{\Psi,j}}_{L^2}^2 + \norm{g_j}_{L^2}^2 - 2\left(f_{\Psi,j}, g_j \circ \phi_{t_j}\right)_{L^2}
  \end{equation*}
  with respect to $t$ is equivalent to maximizing
  \begin{equation*}
    \left(f_{\Psi,j}, g_j \circ \phi_{t_j}\right)_{L^2} = \int_{\R^2} f_{\Psi,j}(y) g_j(y+t_j)\dy = \big(f_{\Psi,j}\star g_j\big)(t_j)
  \end{equation*}
  for all $j\in\{1,\ldots,N\}$. The latter is precisely the method used in the MAL
  algorithm to update the registration.
  
  \paragraph{Remark:} Aside from the reconstruction and registration, a minimizer of $E$
    also yields denoised approximations
    $\mathcal{F}^{-1}\big(\Psi\star_{T_{Z_j}}\Psi\big)$ to the input images for two
    reasons. On the one hand, a minimizer is a least squares approximation according to
    the definition of $E$, which therefore averages the experimental images. On the other
    hand, the simulated real space images $\mathcal{F}^{-1}(\Psi\star_{T_Z}\Psi)$ are
    smooth, which can be seen as follows.
    
    Because of the aperture function and \cref{wcc_l2l2}, the simulated Fourier space
    images $\Psi\star_{T_Z}\Psi$ are compactly supported and bounded for all
    $\Psi\in L^2(A,\C)$ and all $Z\in\R$. By the dominated convergence theorem, the
    partial derivatives are therefore given by
    \begin{equation*}
      \frac{\partial^n}{\partial x^n}\mathcal{F}^{-1}\big(\Psi\star_{T_Z}\Psi\big)(x)
      = \int_{\R^2} \big(\Psi\star_{T_Z}\Psi\big)(y)\frac{\partial^n}{\partial x^n}e^{2\pi\ii x\cdot y}\dy
    \end{equation*}
    for all $n\in\NN^2$.
  
  \medskip
  Using the notation and elementary properties of the weighted cross-correlation from
  \cref{wcc_appendix}, it is easily shown that $\Psi \mapsto E[\Psi, t]$ restricted to
  any one-dimensional affine subspace of its domain $L^2(A, \C)$ is a polynomial of
  degree $4$:
  
  \begin{lemma}\label{subgraphs_1d}
    Let $t\in(\R^2)^N$. For all $\Psi,\Phi\in L^2(A, \C)$ with $\Phi\neq 0$ there are
    coefficients $\big(C_{\Psi,\Phi,t}^j\big)_{j=0,\ldots,4}\in \R^5$ with
    \vspace*{-0.5em}
    \begin{equation*}
      E[\Psi + \alpha\Phi, t] = \sum_{j=0}^4 C_{\Psi,\Phi,t}^j \alpha^j
    \end{equation*}
    for all $\alpha\in\R$. Furthermore, $C_{\Psi,\Phi,t}^4 > 0$.
  \end{lemma}
  \proof Denote the modulated input images by $G_j := \mu_{t_j}\mathcal{F}(g_j)$. Then,
    \begin{align*}
      &E[\Psi + \alpha\Phi, t] = \frac{1}{N}\sum_{j=1}^N\norm{(\Psi+\alpha\Phi)\star_{T_{Z_j}}(\Psi+\alpha\Phi) - G_j}_{L^2}^2 \\
      &= \frac{1}{N}\sum_{j=1}^N\norm{\Psi\star_{T_{Z_j}}\Psi - G_j + \alpha(\Psi\star_{T_{Z_j}}\Phi + \Phi\star_{T_{Z_j}}\Psi) + \alpha^2(\Phi\star_{T_{Z_j}}\Phi)}_{L^2}^2
    \end{align*}
    for all $\alpha\in\R$. Expanding and collecting coefficients with the same power of
    $\alpha$ yields $E[\Psi + \alpha\Phi, t] = \sum_{j=0}^4 C_{\Psi,\Phi,t}^j \alpha^j$
    with
    \begin{align*}
      C_{\Psi,\Phi,t}^0 &= \frac{1}{N}\sum_{j=1}^N \norm{\Psi\star_{T_{Z_j}}\Psi - G_j}_{L^2}^2, \\
      C_{\Psi,\Phi,t}^1 &= \frac{2}{N}\sum_{j=1}^N \Re\Big( \big(\Psi\star_{T_{Z_j}}\Psi - G_j, \Psi\star_{T_{Z_j}}\Phi + \Phi\star_{T_{Z_j}}\Psi\big)_{L^2} \Big), \\
      C_{\Psi,\Phi,t}^2 &= \frac{1}{N}\sum_{j=1}^N \norm{\Psi\star_{T_{Z_j}}\Phi + \Phi\star_{T_{Z_j}}\Psi}_{L^2}^2 + 2\Re\Big( \big(\Psi\star_{T_{Z_j}}\Psi - G_j, \Phi\star_{T_{Z_j}}\Phi\big)_{L^2} \Big), \\
      C_{\Psi,\Phi,t}^3 &= \frac{2}{N}\sum_{j=1}^N \Re\Big( \big(\Psi\star_{T_{Z_j}}\Phi + \Phi\star_{T_{Z_j}}\Psi, \Phi\star_{T_{Z_j}}\Phi\big)_{L^2} \Big), \\
      C_{\Psi,\Phi,t}^4 &= \frac{1}{N}\sum_{j=1}^N \norm{\Phi\star_{T_{Z_j}}\Phi}_{L^2}^2.
    \end{align*}
    Using the symmetry properties in \cref{wcc_symmetry} and the fact that
    $G_j(v) = G_j^*(-v)$ for all $v\in\R^2$, the coefficients of the odd powers of
    $\alpha$ can be simplified to
    \begin{align*}
      C_{\Psi,\Phi,t}^1 &= \frac{4}{N}\sum_{j=1}^N \Re\Big( \big(\Psi\star_{T_{Z_j}}\Psi - G_j, \Psi\star_{T_{Z_j}}\Phi\big)_{L^2} \Big), \\
      C_{\Psi,\Phi,t}^3 &= \frac{4}{N}\sum_{j=1}^N \Re\Big( \big(\Psi\star_{T_{Z_j}}\Phi, \Phi\star_{T_{Z_j}}\Phi\big)_{L^2} \Big).
    \end{align*}
    The coefficient $C_{\Psi,\Phi,t}^4$ is positive by \cref{wcc_injectivity}.\qed
  
  \medskip
  The Gâteaux differentials of $E$ with respect to $\Psi$ can now simply be read off the
  polynomials' coefficients.
  
  \begin{corollary}\label{derivatives_ew}
    Let $t\in(\R^2)^N$ and $\Psi,\Phi\in L^2(A, \C)$ with $\Phi\neq 0$. The Gâteaux
    differentials of $E$ at $(\Psi, t)$ in the direction $(\Phi, 0)$ are
    \begin{align*}
      \big\langle \partial_\Psi E^{(n)}[\Psi,t], \Phi \big\rangle
        := \left.\frac{d}{d\varepsilon} \big\langle \partial_\Psi E^{(n-1)}[\Psi+\varepsilon\Phi,t], \Phi \big\rangle \right|_{\varepsilon=0}
         = \begin{cases} n! C_{\Psi,\Phi,t}^n, & n\in\{1,\ldots,4\}, \\ 0, & n\ge 5. \end{cases}
    \end{align*}
  \end{corollary}
  
  The first order Gâteaux differential of $E$ with respect to the translation $t$ is
  calculated more directly without using difference quotients:
  
  \begin{lemma}\label{derivatives_translation}
    Let $\Psi\in L^2(A,\C)$ and $t,\tilde t\in (\R^2)^N$ with $\tilde t\neq 0$. The
    first order Gâteaux differential of $E$ at $(\Psi, t)$ in the direction
    $(0, \tilde t)$ is
    \begin{equation*}
      \big\langle \partial_t E[\Psi,t], \tilde t\big\rangle
      = -\frac{2}{N}\sum_{j=1}^N\Re\left(\big(\Psi\star_{T_{Z_j}}\Psi - \mu_{t_j}\mathcal{F}(g_j), \nu_{\tilde t_j}\mu_{t_j}\mathcal{F}(g_j)\big)_{L^2}\right),
    \end{equation*}
    where $\nu_{\tilde t_j}(x) = 2\pi i x\cdot \tilde t_j$ for all $x\in\R^2$.
  \end{lemma}
  \proof By the definition of $E$, we have
    \begin{equation*}
      \big\langle \partial_t E[\Psi,t], \tilde t\big\rangle
      = \frac{1}{N}\sum_{j=1}^N\frac{d}{d\varepsilon}\int_{\R^2}\Big|\underbrace{(\Psi\star_{T_{Z_j}}\Psi)(x) - \mu_{t_j+\varepsilon\tilde{t}_j}(x)\mathcal{F}(g_j)(x)}_{=:R_{j,\varepsilon}(x)}\left.\Big|^2\dx\right|_{\varepsilon=0}.
    \end{equation*}
    The partial derivative of the integrand
    is given by
    \begin{equation*}
      \frac{d}{d\varepsilon}\abs{R_{j,\varepsilon}(x)}^2 = \frac{d}{d\varepsilon}\left(R_{j,\varepsilon}(x)R_{j,\varepsilon}^*(x)\right) = 2\Re\left(\left(\frac{d}{d\varepsilon}R_{j,\varepsilon}(x)\right) R_{j,\varepsilon}^*(x) \right),
    \end{equation*}
    where
    \begin{equation*}
      \frac{d}{d\varepsilon}R_{j,\varepsilon}(x) = -\nu_{\tilde t_j}(x)\mu_{t_j+\varepsilon\tilde t_j}(x)\mathcal{F}(g_j)(x).
    \end{equation*}
    This derivative is dominated by
    \begin{equation*}
      \abs{\frac{d}{d\varepsilon}\abs{R_{j,\varepsilon}(x)}^2} \le 2\big|\nu_{\tilde t_j}(x)\big|\abs{\mathcal{F}(g_j)(x)}\left(\abs{\big(\Psi\star_{T_{Z_j}}\Psi\big)(x)} + \abs{\mathcal{F}(g_j)(x)}\right)
    \end{equation*}
    for all $\varepsilon>0$, which is an integrable function by the Cauchy-Schwarz
    inequality, \cref{wcc_lplp} and \cref{functional_definition}. This implies
    \begin{align*}
      \big\langle \partial_t E[\Psi,t], \tilde t\big\rangle
      &= \frac{1}{N}\sum_{j=1}^N\int_{\R^2}\frac{d}{d\varepsilon}\left.\abs{R_{j,\varepsilon}(x)}^2\right|_{\varepsilon=0}\dx \\
      &= -\frac{2}{N}\sum_{j=1}^N\Re\left(\big(\Psi\star_{T_{Z_j}}\Psi - \mu_{t_j}\mathcal{F}(g_j), \nu_{\tilde t_j}\mu_{t_j}\mathcal{F}(g_j)\big)_{L^2}\right). \tag*{\qed}
    \end{align*}
  
  It is apparent that $E$ is not convex with respect to $t$ for arbitrary image series
  $(g_j)_{j=1,\ldots,N}$. The following proposition shows that $E$ is also not convex
  with respect to $\Psi$.
  
  \begin{proposition}\label{convexity}
    Let $t\in (\R^2)^N$. If $g_j\neq 0$ for at least one $j\in\{1,\ldots,N\}$, then
    $\Psi \mapsto E[\Psi, t]$ is not convex.
  \end{proposition}
  \proof Assume that at least one image $g_j$ is nonzero. We consider
    $\Psi\mapsto E[\Psi, t]$ restricted to the lines in $L^2(A, \C)$ that pass through
    the origin. Let $\Phi\in L^2(A, \C)$ with $\Phi\neq 0$. By \cref{subgraphs_1d}, we
    have
    \begin{equation*}
      E[\alpha\Phi, t] = \sum_{j=1}^4 C_{0,\Phi,t}^j \alpha^j
    \end{equation*}
    for all $\alpha\in\R$.
    Since $\Psi=0$, the coefficients of the odd powers of $\alpha$
    are zero and the other coefficients simplify to
    \begin{align*}
      C_{0,\Phi,t}^0 &= \frac{1}{N}\sum_{j=1}^N \norm{g_j}_{L^2}^2, \\
      C_{0,\Phi,t}^2 &= -\frac{2}{N}\sum_{j=1}^N \Re\left(\left( \mu_{t_j}\mathcal{F}(g_j), \Phi\star_{T_{Z_j}}\Phi\right)_{L^2}\right), \\
      C_{0,\Phi,t}^4 &= \frac{1}{N}\sum_{j=1}^N \norm{\Phi\star_{T_{Z_j}}\Phi}_{L^2}^2.
    \end{align*}
    It suffices to show that there is a direction $\Phi\in L^2(A, \C)$ such that the
    biquadratic polynomial
    $E[\alpha\Phi,t] = C_{0,\Phi,t}^4\alpha^4 + C_{0,\Phi,t}^2\alpha^2 + C_{0,\Phi,t}^0$
    is not convex.
    
    The coefficient $C_{0,\Phi,t}^4$ is positive by \cref{subgraphs_1d} and
    $C_{0,\Phi,t}^2$ can be written as
    \begin{equation*}
      C_{0,\Phi,t}^2 = -\frac{2}{N}\sum_{j=1}^N \Re\left(\left(g_j\circ \phi_{t_j}, \mathcal{F}^{-1}\big(\Phi\star_{T_{Z_j}}\Phi\big)\right)_{L^2}\right).
    \end{equation*}
    Assume, we formally extend the weighted cross-correlation and the Fourier transform to the space of tempered distributions.
    Then, we can choose $\Phi = \delta_0$ and get that $\mathcal{F}^{-1}(\Phi\star_{T_{Z_j}}\Phi) = 1$, which corresponds to the simulated TEM
    image of the plane electron wave with no specimen. Then,
    \begin{equation*}
      C_{0,\Phi,t}^2 = -\frac{2}{N}\sum_{j=1}^N \left(g_j\circ\phi_{t_j}, 1\right)_{L^2} = -\frac{2}{N}\sum_{j=1}^N \norm{g_j}_{L^1} < 0,
    \end{equation*}
    which implies that the polynomial $E[\alpha\Phi,t]$ is not convex. The same
    conclusion can be reached without resorting to tempered distributions by considering
    a mollifier sequence for $\Phi$ instead.\qed
  
  \medskip
  Let $t\in(\R^2)^N$. If $g_j = 0$ for all $j\in\{1,\ldots,N\}$, then the functional
  $\Psi\mapsto E[\Psi,t]$ is convex by \cref{auto_wcc_continuousconvex} and $\Psi=0$ is
  the unique global minimizer by \cref{wcc_injectivity}.
  
  \section{Existence of minimizers}\label{existence}
  
  By the direct method of the calculus of variations, minimizers of $E$ exist if $E$ is
  coercive, weakly lower semi-continuous and the domain is reflexive. However, the
  results from \cref{lowpass_filter_appendix} suggest that the functional in its current
  form is not coercive with respect to $\Psi$. For this reason, we add a regularizer to
  $E$ and restrict the admissible set, resulting in
  \begin{equation*}
    E_{\alpha}: L^2(A,\C) \times \overline{B_r(0)}^N \rightarrow \R, \quad (\Psi, t) \mapsto E[\Psi,t] + \alpha\norm{\Psi - \Psi_M}_{L^2}^2 \tag*{$\forall\,\alpha>0$}
  \end{equation*}
  for $\Psi_M\in L^2(A,\C)$ and $r>0$, where $B_r(0)\subseteq\R^2$.
  Instead of regularizing the functional
  with $\alpha\norm{\Psi}_{L^2}^2$, we chose the above generalized nonlinear Tikhonov
  regularization so that the MIMAP functional is included as a special case. The
  regularized functional $E_{\alpha}$ is coercive since
  the translations are bounded, $E$ is non-negative and
  $\alpha\norm{\Psi - \Psi_M}_{L^2}^2 \rightarrow \infty$ as
  $\norm{\Psi}_{L^2}\rightarrow\infty$.
  
  We note that, from a mathematical point of
  view, the additional parameter $r$ must be considered as a workaround to ensure the
  coercivity of the functional in $t$. Unlike usual rigid registration problems, it is not
  easily possible to define an $r>0$ that can serve as an upper bound based on the
  diameter of the image domain. This is because the simulated real space images
  $\mathcal{F}^{-1}(\Psi\star_{T_{Z_j}}\Psi)$ necessarily have unbounded support due to
  the aperture function; therefore it is not obvious whether there exists a minimizing
  sequence of $E$ with bounded translation for all choices of input images
  $g_1,\ldots,g_N$. However, in practice the restriction on the translations is
  justified by the fact that the sample is moving within a bounded domain during image
  acquisition. Therefore, this is a purely mathematical limitation that is 
  irrelevant to the practical application.
  
  \begin{proposition}\label{ew_weak_lower_semicontinuity}
    Let $U\subseteq\R^d$ be bounded, $c\in L^2(U,\C)$ and $w\in W^+(U)$. Then,
    \begin{equation*}
      L^2(U,\C)\times\R^2 \rightarrow \R, \quad (f,t) \mapsto \norm{f\star_w f - \mu_tc}_{L^2}^2
    \end{equation*}
    is weakly lower semi-continuous.
  \end{proposition}
  \proof Split the functional into three parts as follows
    \begin{equation}\label{ew_weak_lower_semicontinuity_splitfunc}
      \norm{f\star_w f - \mu_tc}_{L^2}^2 = \norm{f\star_w f}_{L^2}^2 + \norm{c}_{L^2}^2 - 2\Re\big((f\star_w f, \mu_tc)_{L^2}\big).
    \end{equation}
    The leftmost summand, $f\mapsto \norm{f\star_w f}_{L^2}^2$, is continuous and convex
    by \cref{auto_wcc_continuousconvex} and thus in particular weakly lower
    semi-continuous. The second summand is constant and therefore obviously weakly lower
    semi-continuous.
    
    In order to show that the rightmost
    summand of \cref{ew_weak_lower_semicontinuity_splitfunc} is weakly lower
    semi-continuous, we first consider the special case $w(x,y) = v(x)v^*(y)$ for a
    bounded function $v: U \rightarrow \C$. In the following, we show that the functional
    $G[f,t]:=\big((fv)\star (fv), \mu_tc)_{L^2}$ is weakly continuous on
    $L^2(U,\C)\times\R^2$.
    
    Let $(f, t)\in L^2(U,\C)\times \R^2$ and
    $(f_n, t_n)_{n\in\N}\in(L^2(U,\C)\times \R^2)^\N$ be a weakly convergent sequence
    with $(f_n, t_n) \rightharpoonup (f, t)$. Since
    \begin{align*}
      &\int \int \abs{\big(f_nv\big)^*(y)\big(f_nv)(x+y)(\mu_{t_n}c)^*(x)}\dy\dx \\
      \le\;&\int \norm{f_nv}_{L^2}\norm{(f_nv)(x+\cdot)}_{L^2} \abs{c(x)}\dx \le \norm{f_n}_{L^2}^2\norm{v}_{\infty}^2 \norm{c}_{L^1} < \infty,
    \end{align*}
    we can apply Fubini's theorem to change the integration order so that
    \begin{align*}
      &\big((f_nv)\star (f_nv), \mu_{t_n}c\big)_{L^2} \\
      =\;& \int \big(f_nv\big)^*(y) \int \big(f_nv\big)(x+y) (\mu_{t_n}c)^*(x) \dx \dy \\
      =\;& \int \big(f_nv\big)^*(y) \int f_n(x) \underbrace{v(x)(\mu_{t_n}c)^*(x-y)}_{=:h_{n,y}(x)}\dx \dy = (g_n, f_nv)_{L^2},
    \end{align*}
    where $g_n(y) := \int f_n(x) h_{n,y}(x)\dx$ for all $y\in\R^d$. For
    $h_y(x) := v(x)(\mu_tc)^*(x-y)$, we get
    \begin{equation}
      \abs{h_{n,y}(x) - h_y(x)} \le 2\norm{v}_\infty \abs{c(x-y)} \tag*{$\forall\,x,y\in\R^d$,}
    \end{equation}
    which implies $h_{n,y} \rightarrow h_y$ in $L^2$ for all $y\in\R^d$ by the dominated
    convergence theorem. Since $f_n \rightharpoonup f$, it follows that
    $g_n(y) \rightarrow g(y) := \int f(x)h_y(x)\dx$ for all $y\in\R^d$. Thus, $g_n \rightarrow g$ is a
    pointwise converging sequence.
    
    The support of $(g_n)_{n\in\N}$ is uniformly bounded, since the support of $f_n$ and
    $c$ is contained in $\overline{U}$ for all $n\in\N$. Additionally,
    \begin{equation*}
      \abs{g_n(y)} \le \norm{f_n}_{L^2}\norm{v}_{\infty}\norm{c}_{L^2} \le F \norm{v}_\infty\norm{c}_{L^2}
    \end{equation*}
    holds for all $n\in\N$ and $y\in\R^d$, where
    $F := \sup_{n\in\N} \norm{f_n}_{L^2} < \infty$. Applying the dominated convergence
    theorem once more yields $g_n \rightarrow g$ in $L^2$.
    
    Summing up, $g_n \rightarrow g$ in $L^2$ and $f_nv \rightharpoonup fv$ in $L^2$
    implies $(g_n, f_nv)_{L^2} \rightarrow (g,fv)_{L^2}$ and consequently
    \begin{equation*}
      \big((f_nv)\star (f_nv), \mu_{t_n}c\big)_{L^2} = (g_n, f_nv)_{L^2} \rightarrow (g,fv)_{L^2} = \big((fv)\star (fv), \mu_tc\big)_{L^2}.
    \end{equation*}
    Therefore, $G[f,t]=\big((fv)\star (fv), \mu_tc)_{L^2}$ is weakly
    continuous.
    
    \medskip
    Next, this result is generalized to $H[f,t] := (f\star_w f, \mu_tc)_{L^2}$ using the
    factorization property of the weight $w\in W^+(U)$. By the definition of $W^+(U)$,
    there exist bounded functions $v_{j,N}: U \rightarrow \C$ such that
    $\lim_{N\rightarrow\infty} \norm{w - w_N}_\infty = 0$, where
    \begin{equation*}
      w_N(x,y) = \sum_{j=1}^N v_{j,N}(x)v_{j,N}^*(y) \qquad\forall\,x,y\in U.
    \end{equation*}
    By the previous results,
    \begin{align*}
      H_N: L^2(U,\C)\times\R^2 \rightarrow \C, \quad (f,t) &\mapsto (f\star_{w_N}f, \mu_tc)_{L^2} = \sum_{j=1}^N \big((fv_j)\star(fv_j), \mu_tc\big)_{L^2}
    \end{align*}
    is weakly continuous for all $N\in\N$.
    
    Let $(f, t)\in L^2(U,\C)\times \R^2$ and
    $(f_n, t_n)_{n\in\N}\in(L^2(U,\C)\times \R^2)^\N$ be a weakly convergent sequence
    with $(f_n, t_n) \rightharpoonup (f,t)$. Then,
    \begin{align*}
      &\abs{H[f_n,t_n] - H[f,t]} \\
      =\;&\abs{H[f_n,t_n] - H_N[f_n,t_n] + H_N[f_n,t_n] - H_N[f,t] + H_N[f,t] - H[f,t]} \\
      \le\;&\abs{H[f_n,t_n] - H_N[f_n,t_n]} + \abs{H_N[f_n,t_n] - H_N[f,t]} + \abs{H_N[f,t] - H[f,t]}
    \end{align*}
    holds for all $N\in\N$. Let $C:=\sqrt{\text{Vol}(U)}$ and
    $F := \sup_{n\in\N} \norm{f_n}_{L^2} < \infty$. Then
    $\norm{g}_{L^1} \le C \norm{g}_{L^2}$ for all $g\in L^2(U,\C)$ and the summands in
    the above equation can be estimated by
    \begin{align*}
      \abs{H[f_n,t_n] - H_N[f_n,t_n]} &\le 2\norm{f_n\star_{(w-w_N)}f_n}_{L^2}\norm{c}_{L^2} \\
      \text{(\cref{wcc_l1lp})}&\le 2\norm{f_n}_{L^1}\norm{f_n}_{L^2}\norm{w-w_N}_\infty\norm{c}_{L^2} \\
      &\le 2CF^2\norm{w-w_N}_\infty\norm{c}_{L^2}
    \end{align*}
    and, similarly,
    \begin{equation*}
      \abs{H_N[f,t] - H[f,t]} \le 2CF^2\norm{w-w_N}_\infty\norm{c}_{L^2}.
    \end{equation*}
    Therefore, combined with $\abs{H_N[f_n,t_n] - H_N[f,t]}\to0$ for $n\to\infty$, we get
    \begin{equation*}
      \lim_{n\rightarrow\infty} \abs{H[f_n,t_n] - H[f,t]} \le 4CF^2\norm{w-w_N}_\infty\norm{c}_{L^2}
    \end{equation*}
    for all $N\in\N$. Since  $\lim_{N\rightarrow\infty}\norm{w - w_N}_\infty = 0$ by the
    choice of the sequence $(w_N)_{N\in\N}$, we can conclude
    $\lim_{n\rightarrow\infty}H[f_n,t_n] = H[f,t]$. Thus $H$ is weakly continuous, which
    implies in particular that $(f,t)\mapsto- 2\Re\big((f\star_w f, \mu_tc)_{L^2}\big)$
    is weakly lower semi-continuous. Now the claim follows, since finite sums of weakly
    lower semi-continuous functions are also weakly lower semi-continuous.\qed
  
  \medskip
  Since norms are weakly lower semi-continuous on their respective spaces, the
  regularizer
  \begin{equation*}
    L^2(A,\C) \rightarrow \R, \quad \Psi \mapsto \alpha\norm{\Psi-\Psi_M}_{L^2}^2
  \end{equation*}
  is weakly lower semi-continuous. Overall, $E_{\alpha}$ is weakly lower
  semi-continuous for all $\alpha>0$. Since $L^2(A,\C) \times (\R^2)^N$ is
  a reflexive Banach space and $L^2(A,\C) \times \overline{B_r(0)}{}^N$, as a closed and
  convex subset, is also weakly sequentially closed, the existence of minimizers of
  $E_{\alpha}$ now follows with the direct method:
  \begin{theorem}
    Let $\alpha>0$. There exist $\Psi_*\in L^2(A,\C)$ and $t_*\in \overline{B_r(0)}{}^N$
    with
    \begin{equation*}
      E_{\alpha}[\Psi_*, t_*] = \inf_{\Psi\in L^2(A,\C),\, t\in \overline{B_r(0)}^N} E_{\alpha}[\Psi,t].
    \end{equation*}
  \end{theorem}
  
  \section{Numerical experiment on synthetic data}\label{experiment}
  
  In order to verify that our objective functional $E_{\alpha}$ can indeed be used for
  exit wave reconstruction, we performed a numerical experiment on simulated input data.
  The input images are simulated using the
  same parameters and forward model as in the reconstruction, but on a twice as large
  image area than the area that is used for the reconstruction. The exit wave and images are discretized
  on a cartesian grid with grid width $h=1$ and piecewise constant values for each grid
  cell; the continuous Fourier transform is replaced with the discrete Fourier transform.
  As the algorithm for the numerical minimization of the objective functional we used a
  first order nonlinear Fletcher-Reeves conjugate gradient descent with Armijo step size
  control.
  
  For computational efficiency, we use the TCC $T_Z^\text{Ishizuka}$ given in
  \cref{MAL_tcc} with the focal integration approximation from the MAL algorithm. The
  spatial coherence envelope is then approximated by $E_s(v,w)\approx E_s(v,0)E_s^*(w,0)$,
  whereas the temporal coherence envelope is approximated by a finite sum of factorizable
  terms that originate from a numerical integration of \cref{tcc_mal_et_int}. It is also
  possible to perform the reconstruction without using the focal integration
  approximation, but in this case the TCC can not be written as a finite
  sum of factorizable terms. Hence, it is not possible to utilize the fast Fourier
  transform to speed up the calculation of the simulated images
  $\Psi\star_{T_Z^{\text{Ishizuka}}}\Psi$, which has a significant impact on the
  computation time that is needed to evaluate $E_{\alpha}$ and its derivative.
  
  The exit wave and the focus image series are initially simulated on a grid of
  size $2048\times2048$, which corresponds to an image area of $6.4\times 6.4$ nm$^2$ in
  real space coordinates. Afterwards, the images are cropped to the central section of
  size $1024\times 1024$ pixels, which is then used for the reconstruction. This is done
  in order to reduce the effect of the wrap-around error in the simulated input images,
  which is caused by the periodicity of the discrete fourier transform. Except for this
  additional step, the input images are simulated using the same forward model as in the
  reconstruction. Due to the periodicity of
  our simulated images, the wrap-around
  error only has a small effect and is not accounted for during the reconstruction.
  However, for non-periodic images the wrap around error poses a problem to the numerical
  implementation of the minimization algorithm, which can be seen as a tradeoff with the
  enormous gain of efficiency by using the fast Fourier transform to compute the autocorrelation. In this case, additional
  steps must be taken in order to reduce the effect of the wrap-around error on the reconstructed exit wave.
  See \cite{thust96} for a discussion on how to treat the wrap-around error within the
  context of the MAL algorithm.
  
\begin{figure}
    \centering
    \begin{minipage}[t]{0.35\textwidth}
      \includegraphics[width=\textwidth, keepaspectratio]{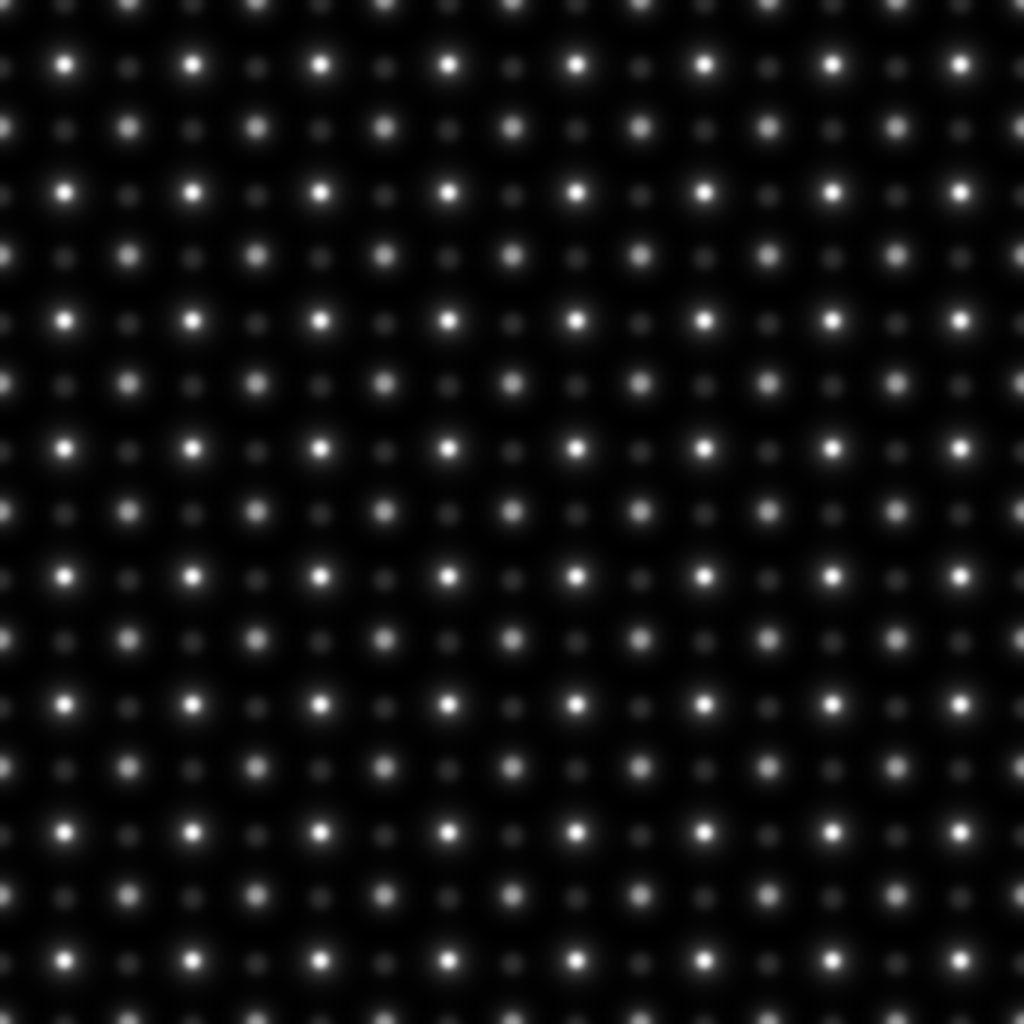}
    \end{minipage}
    \begin{minipage}[t]{0.35\textwidth}
      \includegraphics[width=\textwidth, keepaspectratio]{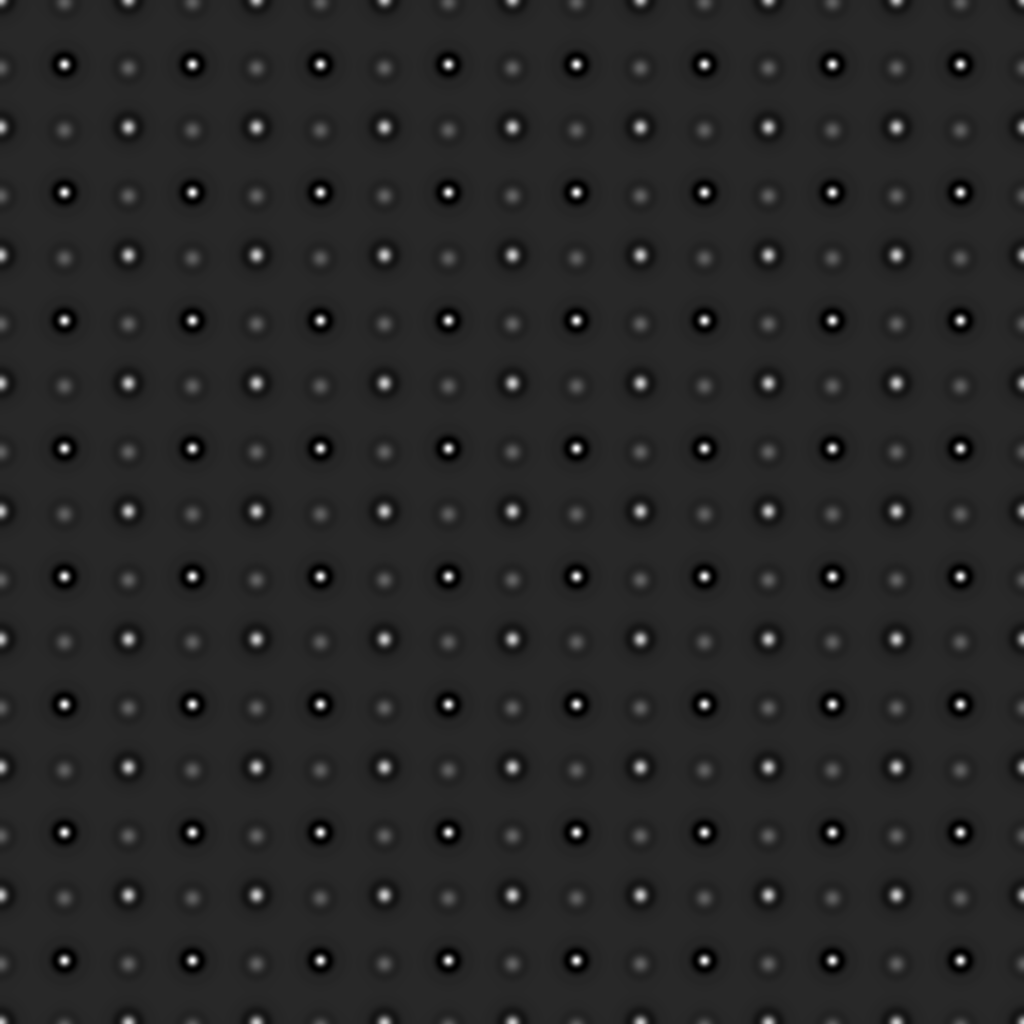}
    \end{minipage}
    \caption{Phase and amplitude of the simulated exit wave. The images depict an area of
    $3.2 \times 3.2$ nm$^2$ and are sampled with $1024\times1024$ pixels.}\label{fig_ew}
  \end{figure}

\begin{figure}
    \centering
    \begin{minipage}[t]{0.4\textwidth}
      \includegraphics[width=\textwidth, keepaspectratio]{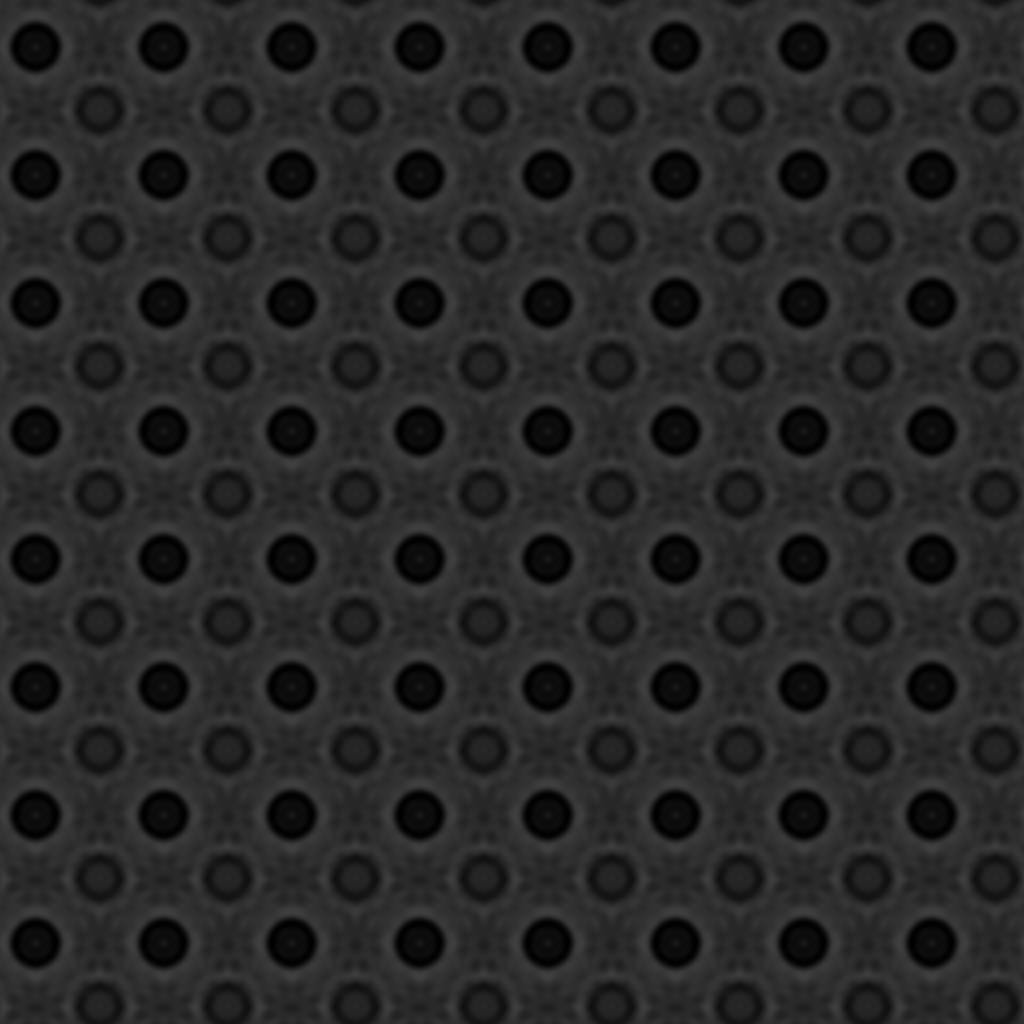}
    \end{minipage}
    \begin{minipage}[t]{0.4\textwidth}
      \includegraphics[width=\textwidth, keepaspectratio]{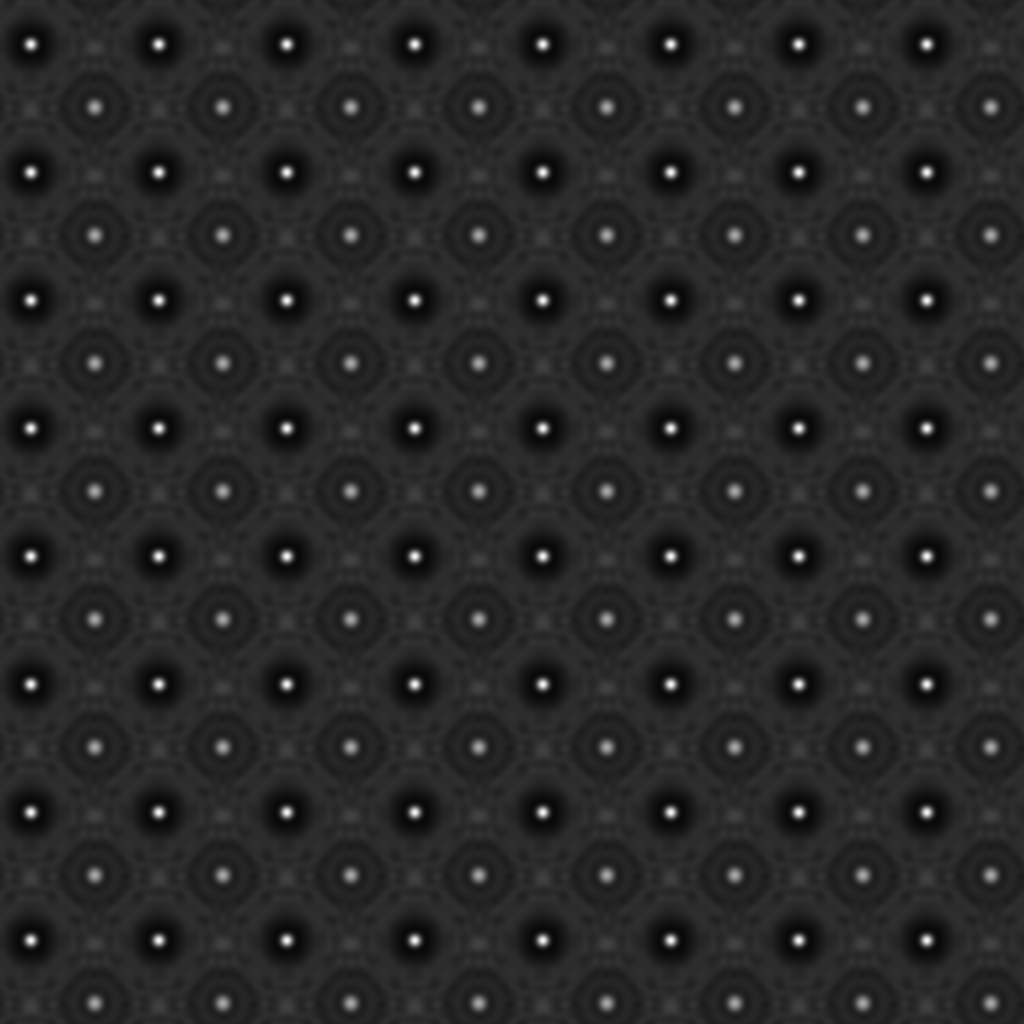}
    \end{minipage}\vspace*{0.2em}
    \begin{minipage}[t]{0.4\textwidth}
      \includegraphics[width=\textwidth, keepaspectratio]{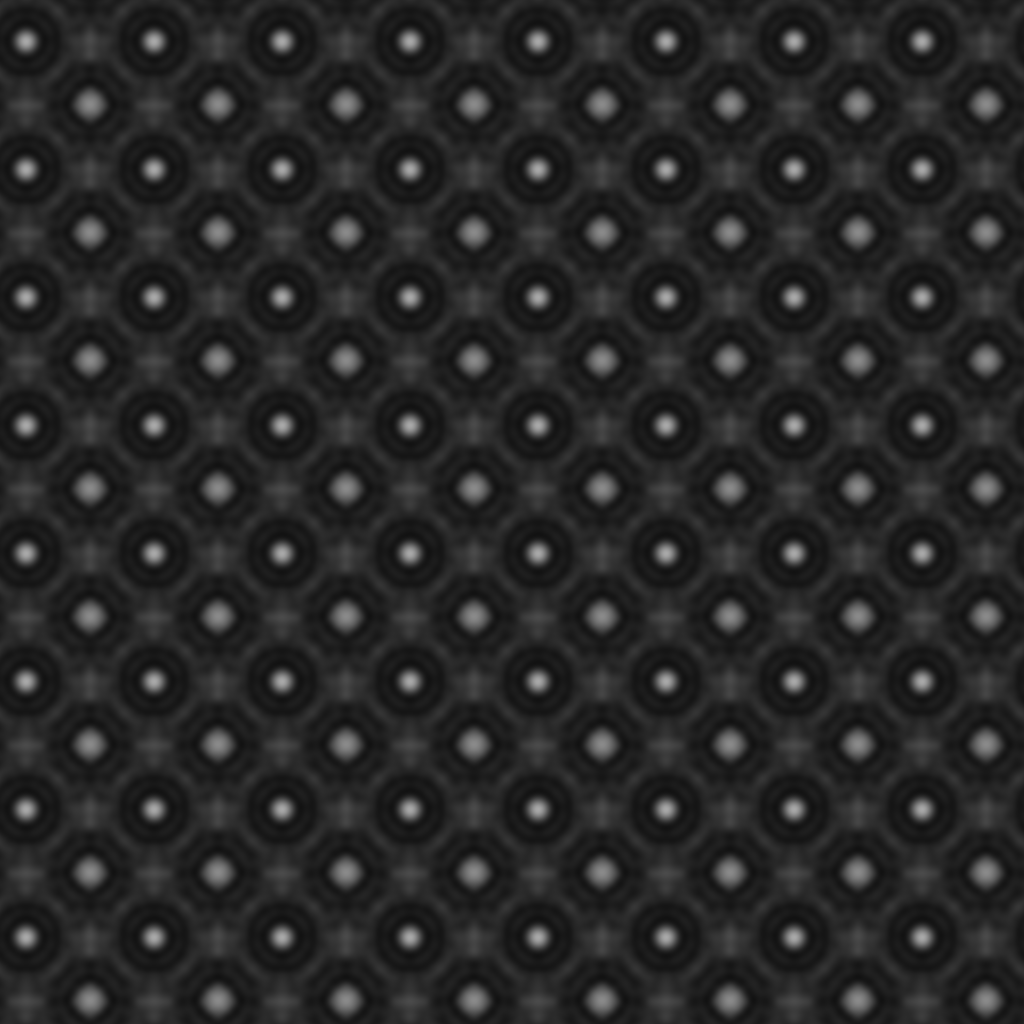}
    \end{minipage}
    \begin{minipage}[t]{0.4\textwidth}
      \includegraphics[width=\textwidth, keepaspectratio]{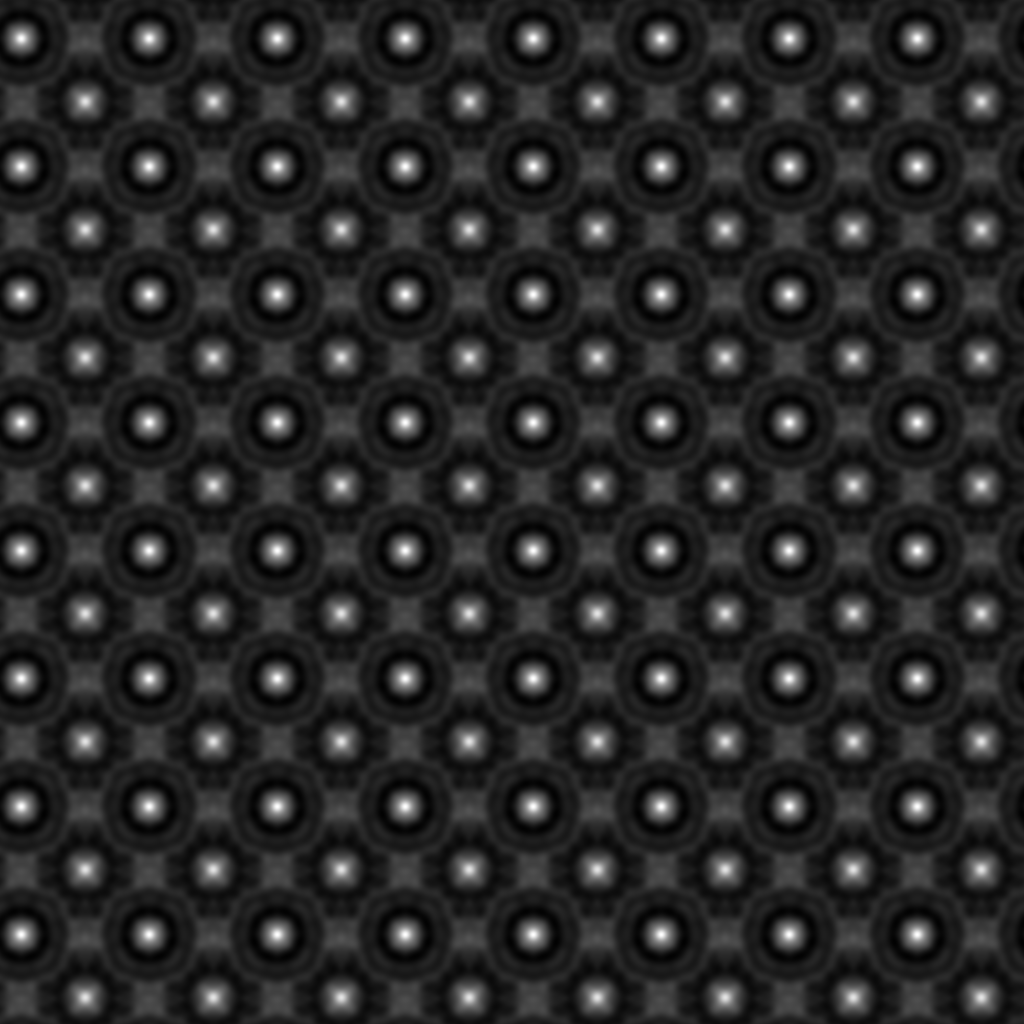}
    \end{minipage}
    \caption{From left to right, top to bottom: Images 6 to 9 of the simulated focus
    series with focus values of $-1$, $0.5$, $2$, and $3.5$ nanometer respectively. The
    image series was simulated with a spherical aberration coefficient of $C_s=-70$nm and
    an objective aperture of $\alpha_\text{max}=125$mrad. The electron wavelength is
    $0.00196875$nm, which corresponds to an accelerating voltage of $300$kV. The images
    depict an area of $3.2\times 3.2$ nm$^2$ and are sampled with $1024\times1024$
    pixels.}\label{fig_simulated_focus_series}
  \end{figure}

  \cref{fig_ew} shows the phase and amplitude of the exit wave that was used to
  simulate the input images for the numerical minimization. The exit wave
  corresponds to a cubic BaTiO$_3$ crystal lattice and was simulated with the multislice
  method \cite{cowley57} using the Dr.\ Probe
  software \cite{barthel18}.
  An excerpt of the simulated focus series is shown in \cref{fig_simulated_focus_series}.
  In total, the series consists of $24$ images sampled with $1024\times 1024$ pixels,
  where the focus values range from $-10$nm to $24.5$nm with a constant focus shift of
  $1.5$nm between successive images in the series. The images were simulated with a
  specimen drift of $0.017$nm between successive images, which corresponds to roughly
  $6$ pixels in the discretized images.
  
  \medskip
  For the practical minimization of the functional $E_{\alpha}$ it is important to
  note two sources of
  non-uniqueness for a minimizer $(\Psi, t)$ in the case $\Psi_M = 0$.
  In this case, multiplying the exit wave with a global phase factor $e^{\ii c}$
  for a $c\in\R$ does not change the energy, i.e. 
  \begin{equation}\label{ambiguity_eq1}
    E_{\alpha}[\Psi e^{ic}, t] = E_{\alpha}[\Psi, t]
  \end{equation}
  holds for all $c\in\R$. However, in practice we observed that the lowest frequency
  of $\Psi$ remained approximately equal to $1$ during the entire minimization,
  which indicates that \cref{ambiguity_eq1} does not affect a derivative based
  minimization method such as the conjugate gradient method. Another source of
  non-uniqueness is the fact that a "global" modulation of
  the exit wave and all images does not change the energy, i.e.
  \begin{equation*}
    E_{\alpha}[\Psi \mu_s, t_s] = E_{\alpha}[\Psi, t]
  \end{equation*}
  holds for all $s\in\R^2$, where
  $t_s = (t_1+s,\, t_2+s,\, ...,\, t_N+s)\in (\R^2)^N$. This ambiguity can easily be
  solved by choosing the first image as a reference and keeping $t_1 = 0$ fixed.
  
  \medskip
  The minimization of our objective functional $E_{\alpha}$ was performed with the
  regularization coefficient set to $\alpha:=10^{-5}$. No a-priori estimate of the
  exit wave is used for the minimization and thus $\Psi_M$ is set to zero. As the initial
  guess for the exit wave we used a constant wave equal to the square root of the mean
  image intensities. The initial guess for the translations was obtained by registering
  successive images in the simulated image series with the cross-correlation. This only
  gives a very rough estimate of the true image shifts, but is sufficiently accurate as
  an initial guess for the minimization. After $183$ steps of the nonlinear conjugate
  gradient method the minimization stops, where we used
  \begin{equation*}
    E_\alpha[\Psi_{k-1},t_{k-1}] - E_\alpha[\Psi_k,t_k] < \varepsilon := 10^{-10}
  \end{equation*}
  as the stopping criterion. Here $(\Psi_k,t_k)$ is the estimate of the exit wave and
  the translation after the $k$-th step.
  
\begin{figure}
    \centering
    \begin{tikzpicture}
      \begin{axis}[width = \textwidth, height = 25em,
                   ymode = log, scaled x ticks = false, tick align = outside, ymajorgrids, xmajorgrids, tick pos = left,
                   xmin = 0, xmax = 183,
                   ymin = 1e-8, ymax = 1e-1,
                   xlabel=number of iterations]
        \addplot[mark=none, color=orange] table [x expr=\thisrowno{0}, y expr=\thisrowno{1}*1024^2] {energy.dat};
        \addlegendentry{$E_{\alpha}[\Psi_k,t_k]$}
        \addplot[mark=none, color=blue] table [x expr=\thisrowno{0}, y expr=\thisrowno{1}*1024^2] {dataterm.dat};
        \addlegendentry{$E_{\alpha}[\Psi_k,t_k] - \alpha\norm{\Psi_k}_{L^2}^2$}
        \addplot[mark=none, color=cyan] table [x expr=\thisrowno{0}, y expr=\thisrowno{1}*1024^2] {regularizer.dat};
        \addlegendentry{$\alpha\norm{\Psi_k}_{L^2}^2$}
      \end{axis}
    \end{tikzpicture}
    \caption{Logarithmic plot of the energy and its components on the $y$-axis and the number of
    iterations of the nonlinear conjugate gradient method on the $x$-axis.}\label{fig_energy}
  \end{figure}
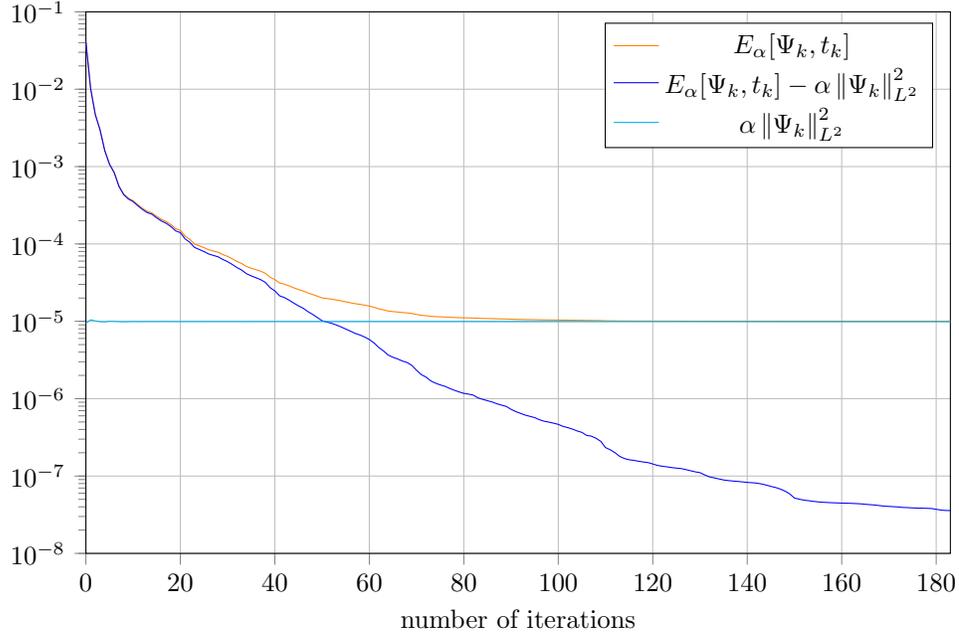

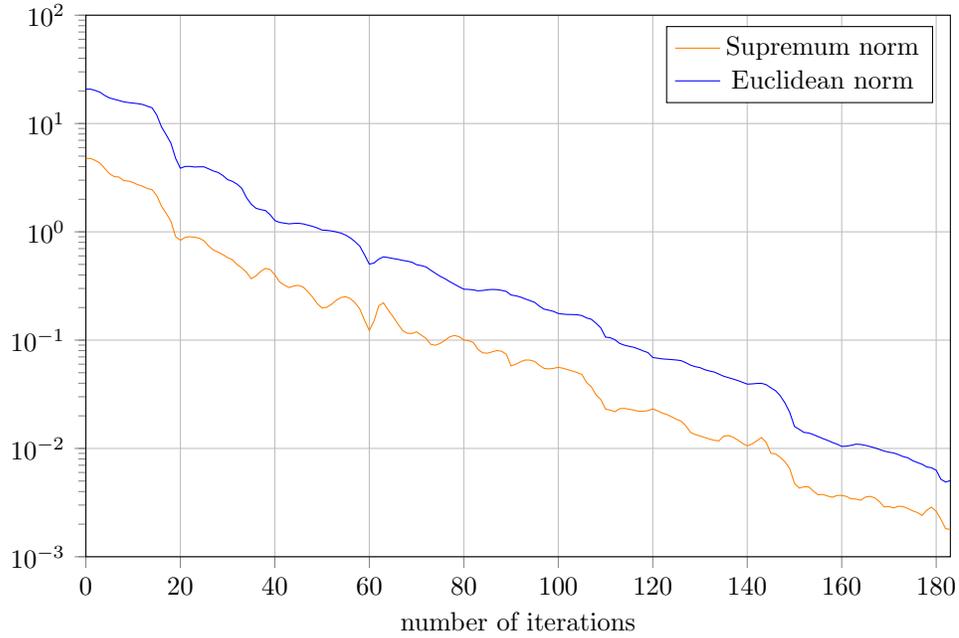
\begin{figure}
    \centering
    \begin{tikzpicture}
      \begin{axis}[width = \textwidth, height = 25em,
                   ymode = log, scaled x ticks = false, tick align = outside, ymajorgrids, xmajorgrids, tick pos = left,
                   xmin = 0, xmax = 183,
                   ymin = 1e-3, ymax = 1e2,
                   xlabel=number of iterations]
        \addplot[mark=none, color=orange] table [x expr=\thisrowno{0}, y expr=\thisrowno{1}*1024/3.2] {translation_sup.dat};
        \addlegendentry{Supremum norm}
        \addplot[mark=none, color=blue] table [x expr=\thisrowno{0}, y expr=\thisrowno{1}*1024/3.2] {translation_euc.dat};
        \addlegendentry{Euclidean norm}
      \end{axis}
    \end{tikzpicture}
    \caption{Distance of the estimated translation to the correct translation with
    respect to the supremum norm (orange graph) and the euclidean norm (blue graph). The
    number of iterations of the nonlinear conjugate gradient method is given on the
    $x$-axis and the distance in pixels is given on the $y$-axis.}\label{fig_translation}
  \end{figure}

\begin{figure}
    \centering
    \begin{tikzpicture}
      \begin{axis}[width = \textwidth, height = 25em,
                   ymode = log, scaled x ticks = false, tick align = outside, ymajorgrids, xmajorgrids, tick pos = left,
                   xmin = 0, xmax = 183,
                   ymin = 1e-2, ymax = 1e3,
                   xlabel=number of iterations]
        \addplot[mark=none, color=orange] table{ew_realspace_sup_residual.dat};
        \addlegendentry{Supremum norm}
        \addplot[mark=none, color=blue] table{ew_realspace_euc_residual.dat};
        \addlegendentry{Euclidean norm}
      \end{axis}
    \end{tikzpicture}
    \caption{Distance of the estimated exit wave in real space coordinates to the central
    section of the exit wave that was used for the input data simulation with respect to
    the supremum norm (orange graph) and the euclidean norm (blue graph, not normalized by
    the number of pixels).}\label{fig_ewplot}
  \end{figure}
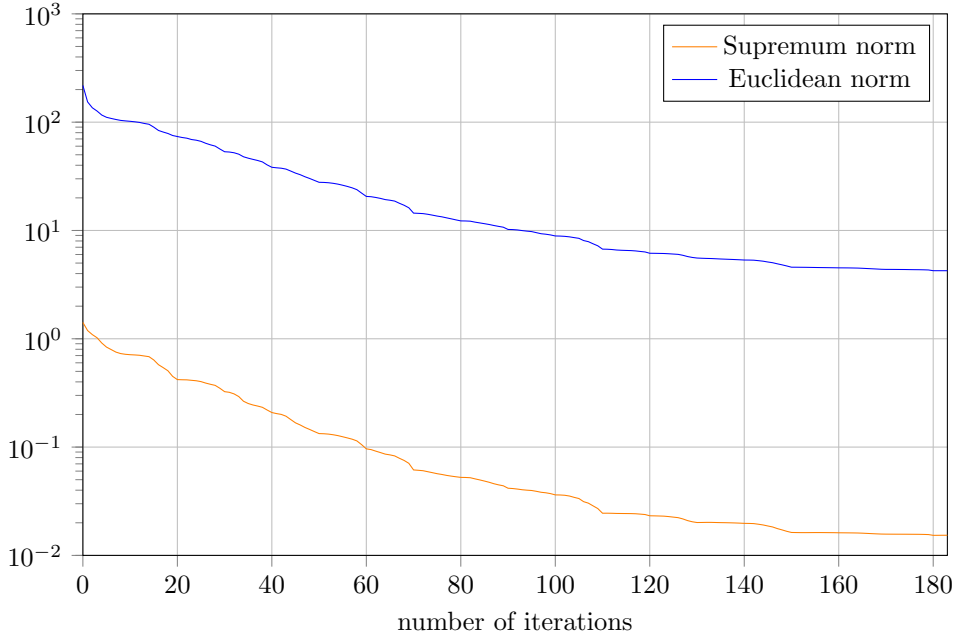

  In \cref{fig_energy} it can be seen
  that the data term is minimized successfully, while the value of the regularizer
  essentially stays constant. Furthermore, \cref{fig_translation} shows that the
  minimization also yielded good approximations to the correct translations and
  \cref{fig_ewplot} shows that the reconstructed exit wave is an approximation to
  the central section of the exit wave that was used for the input data simulation.
  Neither the phase and the amplitude of the reconstructed exit wave nor simulated
  TEM images based on the reconstructed exit wave are shown here, since they are visually
  identical to the images in \cref{fig_ew,fig_simulated_focus_series}. The pixel value
  ranges are $[0.0693, 1.815]$ (reconstructed), $[0.0741, 1.824]$ (correct) of the phase
  and $[0.8747, 1.637]$ (reconstructed), $[0.8746, 1.645]$ (correct) of the amplitude of
  the exit wave.
           
  \section{Conclusions and outlook}\label{outlook}
  
  We have shown that the transmission cross-coefficient $T_Z$ satisfies the factorization
  property, which, together with the concept of weighted cross-correlations, formed the
  basis for all further results. We then
  proposed a novel functional for joint exit wave reconstruction and image registration
  and derived expressions for the first and higher order Gâteaux differentials. We have
  shown that our objective functional is not convex with respect to the exit wave except
  in a trivial case, which also applies to the MAL functional. Using the direct method, we
  continued to show that minimizers of the Tikhonov-regularized version of our objective
  functional exist. Finally, the applicability of our approach was demonstrated with a
  numerical experiment on simulated input data.
  
  In its current form, the problem of exit wave reconstruction using our objective
  functional is clearly not well-posed, since a minimizer is not unique. Therefore, it
  would be interesting to investigate additional assumptions or restrictions on our
  objective functional that are sufficient to ensure a unique minimizer. Similarly, it is
  not clear if the reconstructed exit wave depends continuously on the input data.
  
  In practice, the focus step size between successive images in a focus series is very
  well known, but the available estimate of the actual focus value from the microscope
  settings is quite inaccurate. Therefore, an extension of our model that includes the
  base focus value as an unknown appears to be especially useful for the application of
  our method to real experimental data obtained with a TEM.

  Extending the registration method to a piecewise rigid registration would allow for a
  reconstruction of the exit wave even if the specimen consists of multiple parts that
  move in different directions.
  
  From a theoretical point of view it would be desirable to have a more complete
  result regarding the coercivity, since we have only shown that a particular variant of
  our objective functional is not coercive.
  
  \section*{Acknowledgements}
  The authors were funded in part by the Excellence Initiative of the German Federal and State Governments through grant GSC 111.

  \appendix
  \section{Weighted cross-correlation}\label{wcc_appendix}
  
  This section is a collection of various results regarding a generalization of the
  ordinary cross-correlation. The proofs are largely a direct generalization of the
  corresponding proofs for the ordinary cross-correlation.
  
  \begin{definition}\label{wcc_def}
    Let $f,g: \R^d \rightarrow \C$ and $w: \R^d \times \R^d \rightarrow \C$. The weighted
    cross-correlation $f \star_w g$ is defined as
    \begin{equation*}
      \big(f \star_w g\big)(x) := \int_{\R^d} f^*(y)g(x+y)w(x+y,y)\dy
    \end{equation*}
    for all $x\in\R^d$ such that the integral is well defined and finite.
  \end{definition}
  
  This definition is extended to functions $f: M_1 \rightarrow \C$ and
  $g: M_2 \rightarrow \C$ for subsets $M_1, M_2 \subseteq \R^d$ by setting $f(x) := 0$
  for $x\in\R^d\backslash M_1$ and $g(x) := 0$ for $x\in\R^d\backslash M_2$. It is
  similarly extended to weight functions $w: M_3 \times M_4 \rightarrow\C$ for subsets
  $M_3, M_4\subseteq\R^d$.
  
  Two vector spaces of weight functions are particularly useful here, so they are given a
  special notation. For $U\subseteq\R^d$, the space of measurable and bounded weight
  functions is defined as
  \begin{equation*}
    W(U) := \{ w:U\times U \rightarrow \C \mid w \text{ measurable and bounded}\}.
  \end{equation*}
  The subspace $W^+(U)$ of $W(U)$, where every function can be approximated by a
  particular kind of sequence of factorizable functions, is defined as
  \begin{align*}
    W^+(U) := \Big\{ w\in &W(U) \;\Big\vert\; \forall\,N\in\N\;\forall\,j\in\{1,\ldots,N\}\;\exists\, w_N \in W(U)\\
                 &\;\;\exists\,v_{j,N}: U \rightarrow \C \text{ measurable and bounded}: \\
                 &\left.\lim_{N\rightarrow\infty} \norm{w - w_N}_\infty = 0\;\wedge\; w_N(x,y) = \sum\nolimits_{j=1}^N v_{j,N}(x)v_{j,N}^*(y)\right\}.
  \end{align*}
  
  In the following, we will consider the weighted cross-correlation for Lebesgue
  functions. This makes sense for $f$ and $g$, as the integration in \cref{wcc_def} is
  carried out over the entire domain of $f$ and $g$. However, it does not make sense to
  consider the weight as a Lebesgue function, as the integration is only carried out over
  a subset of the weight's domain of measure zero (namely, the diagonals
  $\{(x+y,y) \mid y\in\R^d\}$). It should be possible to extend the definition to
  weights $\widetilde w\in L^\infty(U\times U, \C)$ by
  treating $f\star_{\widetilde w} g$ itself as a Lebesgue function, but this approach does not have
  any benefits for the theory developed in
  \cref{forward_model,objective_functional,existence} and therefore is not pursued.
  Nevertheless, we occasionally
  identify a weight $w\in W(U)$ with the equivalence class
  $[w]\in L^\infty(U\times U, \C)$ in order to be able to apply Hölder's inequality.
  
  \begin{lemma}\label{wcc_l2l2}
    If $f,g\in L^2(\R^d, \C)$ and $w\in W(\R^d)$, then $(f\star_wg)(x)\in\C$ for all
    $x\in\R^d$ and $\norm{f\star_wg}_\infty\le \norm{w}_\infty\norm{f}_{L^2}\norm{g}_{L^2}$.
  \end{lemma}
  \proof This follows from the Cauchy-Schwarz inequality.\qed
  
  \begin{lemma}\label{wcc_l1lp}
    Let $p\in[1,\infty)$. 
    \hspace*{-0.1em}If $f\in L^1(\R^d, \C)$, $g\in L^p(\R^d, \C)$ and $w\in W(\R^d)$,
    then $f \star_w g \in L^p(\R^d, \C)$ and
    $\norm{f\star_wg}_{L^p} \le \norm{w}_{\infty} \norm{f}_{L^1} \norm{g}_{L^p}$.
  \end{lemma}
  \proof This can be shown with a straightforward application of Minkowski's
    integral inequality.\qed
  
  \medskip
  An immediate consequence of \cref{wcc_l1lp} is that the space $L^1(\R^d, \C)$ is closed
  under weighted cross-correlation. If the support of the weight is bounded, then the
  same is true for every $L^p$-space with $p\in[1,\infty)$. 
  
  \begin{corollary}\label{wcc_lplp}
    If $p\in[1,\infty)$ 
    and $w\in W(U)$ for a bounded subset $U\subseteq\R^d$, then
    $f\star_w g \in L^p(\R^d, \C)$ for all $f,g\in L^p(\R^d,\C)$.
  \end{corollary}
  \proof Since the support of $w$ is bounded, there is an $r>0$ such that
    \begin{equation*}
      \mathrm{supp}(w)\subseteq B_r(0) \subseteq \R^d \times \R^d.
    \end{equation*}
    Let $A := B_r(0) \subseteq \R^d$. Then $w(x,y) = w(x,y) \chi_A(x)\chi_A^*(y)$ for all
    $x,y\in\R^d$ and it follows that
    \begin{equation*}
      f\star_w g = (f\chi_A) \star_w (g\chi_A).
    \end{equation*}
    Since $f\chi_A\in L^p(\R^d, \C)$ has bounded support, we get
    $f\chi_A\in L^1(\R^d, \C)$ from Hölder's inequality and consequently
    $f\star_w g \in L^p(\R^d, \C)$ by the previous lemma. \qed
  
  \medskip
  The next lemma lists two symmetry properties that are useful to simplify calculations
  with weighted cross-correlations.
  
  \begin{lemma}\label{wcc_symmetry}
    Let $f,g\in L^2(\R^d, \C)$ and $w\in W(U)$ for a bounded subset $U\subseteq\R^d$ with
    $w^*(x,y) = w(y,x)$ for all $x,y\in\R^d$. Then
    \begin{enumerate}
      \item $\big(f\star_wg\big)(x) = \big(g\star_wf\big)^*(-x)$ for all $x\in\R^d$.
      \item $(h, f\star_wg)_{L^2} = (g\star_wf, h)_{L^2}$ for all $h\in L^2(\R^d, \C)$
            with $h(x) = h^*(-x)$ for all $x\in\R^d$.
    \end{enumerate}
  \end{lemma}
  \proof Both statements follow directly from the definition of the weighted
    cross-correlation. \qed
  
  \medskip
  The first symmetry property can be used to derive the following results on the weighted
  autocorrelation of $L^2$-functions.
    
  \begin{lemma}\label{auto_wcc_nonnegative}
    Let $f\in L^2(\R^d,\C)$ and $w\in W^+(U)$ for a bounded subset $U\subseteq\R^d$ with
    $w^*(x,y) = w(y,x)$ for all $x,y\in\R^d$. Then $\mathcal{F}^{-1}(f\star_wf)$ is
    real-valued and nonnegative.
  \end{lemma}
  \proof By the first symmetry property of \cref{wcc_symmetry}, we have
    \begin{equation*}
      \big(f\star_wf\big)(x) = \big(f\star_wf\big)^*(-x)
    \end{equation*}
    for all $x\in\R^d$, which implies that $\mathcal{F}^{-1}(f\star_wf)$ is real-valued.
    
    Let $(w_N)_{N\in\N} \in W(U)^\N$ and $v_{j,N}: U \rightarrow \C$ measurable and
    bounded such that $\lim_{N\rightarrow\infty}\norm{w - w_N}_\infty = 0$ and
    $w_N(x,y) = \sum_{j=1}^N v_{j,N}(x)v_{j,N}^*(y)$ for all $N\in\N$. By the same
    reasoning as in \cref{wcc_lplp} it is possible to assume without loss of generality
    that $\supp(f)$ is bounded. This implies $f\in L^1(\R^d, \C)$ and shows that
    $f\star_{w_N} f\in L^1(\R^d, \C)$ for all $N\in\N$ by \cref{wcc_l1lp}. Therefore,
    $\mathcal{F}^{-1}(f \star_{w_N} f)$ is continuous and the inequality
    \begin{equation}\label{auto_wcc_nonnegative_eq}
      \mathcal{F}^{-1}(f \star_{w_N} f) = \sum_{j=1}^N \abs{\mathcal{F}^{-1}(fv_{j,N})}^2 \ge 0
    \end{equation}
    holds for all $N\in\N$ by the linearity of the Fourier transform and the convolution
    theorem. Since the Fourier transform is unitary, it follows that
    \begin{align*}
      \norm{\mathcal{F}^{-1}(f\star_wf) - \mathcal{F}^{-1}(f\star_{w_N}f)}_{L^2}
      &= \norm{f\star_{(w-w_N)}f}_{L^2} \\
      &\le \norm{w-w_N}_{\infty}\norm{f\star f}_{L^2} \underset{N\rightarrow\infty}{\longrightarrow} 0,
    \end{align*}
    where $\norm{f\star f}_{L^2}<\infty$ since $\supp(f)$ is bounded. Thus
    $\mathcal{F}^{-1}(f\star_{w_N}f) \rightarrow \mathcal{F}^{-1}(f\star_wf)$ in $L^2$
    and consequently there exists a subsequence of
    $\big(\mathcal{F}^{-1}(f\star_{w_N}f)\big)_{N\in\N}$ converging pointwise to
    $\mathcal{F}^{-1}(f\star_wf)$. But all elements of
    $\big(\mathcal{F}^{-1}(f\star_{w_N}f)\big)_{N\in\N}$ are nonnegative by
    \cref{auto_wcc_nonnegative_eq}, which shows that $\mathcal{F}^{-1}(f\star_wf)$ is
    nonnegative as well.\qed
  
  \medskip
  \begin{corollary}\label{auto_wcc_continuousconvex}
    Let $U\subseteq\R^d$ be bounded and $w\in W^+(U)$ such that $w^*(x,y) = w(y,x)$ for
    all $x,y\in\R^d$. Then the functional
    \begin{equation*}
      F: L^2(\R^d,\C) \rightarrow \R, \quad f \mapsto \norm{f\star_wf}_{L^2}^2
    \end{equation*}
    is Fréchet differentiable, continuous and convex.
  \end{corollary}
  \proof Analogously to \cref{derivatives_ew}, we get that
    \begin{equation*}
      F'[f]: L^2(\R^d, \C) \rightarrow \R, \quad g \mapsto 2\Re\big((f\star_wf, f\star_wg+g\star_wf)_{L^2}\big)
    \end{equation*}
    is the first order Gâteaux differential of
    $F$ at $f\in L^2(\R^d,\C)$ in the direction $g\in L^2(\R^d,\C)$ and
    \begin{align*}
      \abs{F'[f](g)} &\le 2\norm{f\star_wf}_{L^2}\norm{f\star_wg+g\star_wf}_{L^2} \le 4\norm{f\star_wf}_{L^2}\norm{f\star_wg}_{L^2} \\
      &\le 4 \norm{f\star_wf}_{L^2}\norm{w}_\infty\norm{f}_{L^1}\norm{g}_{L^2}
    \end{align*}
    holds for all $g\in L^2(U,\C)$ by \cref{wcc_symmetry} (i) and \cref{wcc_l1lp}.
    Therefore, $F'[f]$ is a bounded linear operator that additionally satisfies
    \begin{equation*}
      \lim_{g\rightarrow 0} \frac{\abs{F[f+g]-F[f]-F'[f](g)}}{\norm{g}_{L^2}} = 0,
    \end{equation*}
    which shows that $F$ is Fréchet differentiable and particularly continuous.
    
    Analogously to \cref{derivatives_ew}, the second order Gâteaux differential of $F$ at
    a position $f\in L^2(U,\C)$ in the direction $g\in L^2(U,\C)$ is
    \begin{align*}
      \langle F''[f], g\rangle &= \norm{f\star_wg+g\star_wf}_{L^2}^2 + 2\Re\big((f\star_wf, g\star_wg)_{L^2}\big) \\
      &= \norm{f\star_wg+g\star_wf}_{L^2}^2 + 2\Re\big((\mathcal{F}^{-1}(f\star_wf), \mathcal{F}^{-1}(g\star_wg))_{L^2}\big).
    \end{align*}
    By \cref{auto_wcc_nonnegative} we have $\langle F''[f], g\rangle \ge 0$ for all
    $f,g\in L^2(U, \C)$, which implies that $F$ is convex.\qed
  
  \medskip
  Similar to the ordinary cross-correlation, the weighted cross-correlation of $L^p$- and
  $L^q$-functions is continuous for suitable weight functions. Note that the constant
  weight function $w(x,y) = 1$ for all $x,y\in\R^d$ satisfies all of the conditions in
  the following Lemma, which therefore is a direct generalization of the corresponding
  statement for the ordinary cross-correlation.
  
  \begin{lemma}\label{wcc_continuity}
    Let $U\subseteq\R^d$ be an open set such that the measure of its boundary
    $\partial U$ is zero and $p,q\in(1,\infty)$ with
    $\frac{1}{p}+\frac{1}{q} = 1$. Choose $f\in L^p(U, \C)$, $g\in L^q(U, \C)$
    and $w\in W(U)$. If $w$ is continuous on $U\times U$, then
    $f \star_w g\in C(\R^d,\C)$.
  \end{lemma}
  \proof Fix an arbitrary $x\in\R^d$ and let $\varepsilon>0$. We show that $f\star_wg$ is
    continuous at $x$. Since $C_c(\R^d,\C)$ is dense in $L^q(\R^d,\C)$, there is a
    function $v\in C_c(\R^d,\C)$ with $\norm{g-v}_{L^q} < \varepsilon$. Define
    \begin{align*}
      &G_z: \R^d \rightarrow \C, \quad y \mapsto g(z+y)w(z+y,y), \\
      &V_z: \R^d \rightarrow \C, \quad y \mapsto v(z+y)w(z+y,y)
    \end{align*}
    for all $z\in\R^d$. By the Hölder and Minkowski inequalities,
    \begin{align*}
      &\abs{\big(f\star_w g\big)(x+h) - \big(f\star_w g\big)(x)} \\
      \le\;&\int \abs{ f^*(y)g(x+h+y)w(x+h+y,y) - f^*(y)g(x+y)w(x+y,y)}\dy \\
      \le\;&\norm{f}_{L^p} \norm{ G_{x+h} - G_x }_{L^q} \\
      \le\;&\norm{f}_{L^p} \big(\norm{ G_{x+h} - V_{x+h} }_{L^q} + \norm{ V_{x+h} - V_x }_{L^q} + \norm{V_x - G_x}_{L^q}\big)
    \end{align*}
    holds for all $h\in\R^d$. An upper bound for the first and third summand is given by
    \begin{equation}
      \norm{G_z - V_z}_{L^q} \le \norm{g-v}_{L^q} \norm{w}_{\infty} < \varepsilon \norm{w}_{\infty} \tag*{$\forall\,z\in\R^d$.}
    \end{equation}
    By the continuity of $v$ and the continuity of $w$ on $U\times U$, we have
    \begin{align*}
      \lim_{h\rightarrow 0} V_{x+h}(y) &= \lim_{h\rightarrow 0} v(x+h+y)w(x+h+y,y) \\
      &= \begin{cases}
           v(x+y)w(x+y,y), & \text{if $y\in U \;\;\wedge\;\; x+y\in U$}, \\
           0,              & \text{if $y\in U \;\;\wedge\;\; x+y\in\overline{U}^c$}, \\
           0,              & \text{if $y\in U^c$},
         \end{cases} \\
      &= \begin{cases}
           V_x(y), & y\in U \;\;\wedge\;\; x+y\in U, \\
           V_x(y), & y\in U \;\;\wedge\;\; x+y\in\overline{U}^c, \\
           V_x(y), & y\in U^c.
         \end{cases}
    \end{align*}
    This implies $\lim_{h\rightarrow 0} V_{x+h}(y) = V_x(y)$ for almost all $y\in\R^d$,
    since $\partial U$ has zero measure and thus also the set
    $\{y\in U \mid x+y\in\partial U\}$. Furthermore, since
    $\abs{V_{x+h}(y)} \le \norm{w}_{\infty}\norm{v}_{\infty} < \infty$ for all $y\in\R^d$
    and
    \begin{equation*}
      \supp(V_{x+h}) \subseteq \overline{B_1\big(\supp(v) - x\big)} =: D \qquad\qquad\forall\,h\in B_1(0),
    \end{equation*}
    it follows that $\norm{w}_{\infty}\norm{v}_{\infty} \chi_D$ is an integrable
    function that dominates $V_{x+h}$ for all $h\in B_1(0)$. Now the dominated
    convergence theorem implies
    \begin{equation*}
      \lim_{h\rightarrow 0} \norm{ V_{x+h} - V_x }_{L^q} = 0
    \end{equation*}
    and we conclude $\lim\limits_{h\rightarrow 0}\abs{\big(f\star_w g\big)(x+h) - \big(f\star_w g\big)(x)} = 0$.
    \qed
  
  \medskip
  An application of the continuity of the weighted cross-correlation is the following
  characterization of $L^2$-functions, whose weighted autocorrelation is zero.
  
  \begin{lemma}\label{wcc_injectivity}
    Let $U\subseteq\R^d$ be an open set such that the measure of its boundary
    $\partial U$ is zero. Choose $f\in L^2(U,\C)$ and a continuous weight
    $w\in W(U)$. If there is a $c>0$ with $w(y,y) \ge c$ for almost all
    $y\in U$, then
    \begin{equation*}
      \norm{f\star_wf}_{L^2} = 0 \qquad\Longleftrightarrow\qquad \norm{f}_{L^2} = 0.
    \end{equation*}
  \end{lemma}
  \proof "$\Leftarrow$": $\norm{f}_{L^2} = 0$ implies $f = 0$ almost everywhere and thus
    $\norm{f\star_wf}_{L^2} = 0$.
    "$\Rightarrow$": If $\norm{f\star_wf}_{L^2} = 0$, then $f\star_wf = 0$ almost
    everywhere. However, the function $f\star_wf$ is continuous by \cref{wcc_continuity},
    which implies $f\star_wf = 0$. In particular,
    \begin{equation*}
      0 = \big(f\star_wf\big)(0) = \int_U\abs{f(y)}^2w(y,y)\dy \ge c\norm{f}_{L^2}^2 \ge 0
    \end{equation*}
    and thus $\norm{f}_{L^2} = 0$.\qed
  
  \section{Low-pass filter approximation}\label{lowpass_filter_appendix}
  
  It is easy to see that the functional $f \mapsto \norm{f\star f}_{L^2}$  with
  $f\in L^2(\R^d, \C)$ is not coercive. If we define
  \begin{equation*}
    g_{\delta}: \R^d \rightarrow \R, \quad x \mapsto \begin{cases}
      \prod_{i=1}^d \frac{1}{\sqrt{x_i}}, & \text{if } 1\le x_i\le\delta \;\,\forall\,i\in\{1,\ldots,d\}, \\
      0, & \text{otherwise}
    \end{cases}
  \end{equation*}
  for all $\delta>1$, then $\lim_{\delta\rightarrow\infty}\norm{g_\delta}_{L^2} = \infty$
  and $\lim_{\delta\rightarrow\infty}\norm{g_\delta^2}_{L^2} = 1$. Therefore, it follows
  that
  \begin{equation*}
    \lim_{\delta\rightarrow\infty}\norm{f_\delta}_{L^2} = \infty\qquad\text{and}\qquad\lim_{\delta\rightarrow\infty}\norm{f_\delta\star f_\delta}_{L^2} = \lim_{\delta\rightarrow\infty}\norm{(\mathcal{F}^{-1}f_\delta)^2}_{L^2} = 1
  \end{equation*}
  for $f_\delta = \mathcal{F}g_\delta$.
  
  The purpose of this section is to generalize this result to
  \begin{equation*}
    f \mapsto \norm{f\star f}_{L^2}, \qquad f\in L^2(U, \C)
  \end{equation*}
  for any open set $U\subseteq\R^d$ with non-empty interior. The central tool to this end
  is \cref{convolution_approximation}, which provides a method to estimate how well a
  function $g:\R\rightarrow\R$ is approximated by $g * \sinc$. Here, $\sinc$ is the
  normalized cardinal sine function defined as
  \begin{equation*}
    \sinc(x) = \begin{cases}
                     \frac{\sin(\pi x)}{\pi x}, & \text{if $x\neq 0$,} \\
                     1,                         & \text{if $x=0$.}
                 \end{cases}
  \end{equation*}
  We have $\mathcal{F}(\sinc) = \chi_{[-\frac{1}{2},\frac{1}{2}]}$ so that the
  convolution of $g$ with $\sinc$ corresponds to applying an ideal low-pass filter with
  cutoff frequency $\frac{1}{2}$ to $g$.
  
  In order to compute bounds for real-valued functions, we make use of interval
  arithmetic. The sum and the difference of two intervals $I_1 := [a,b]$ and
  $I_2 := [c,d]$ are
  \begin{align*}
    I_1 + I_2 &:= \{x+y\mid x\in I_1,\,y\in I_2\} = [a+c,b+d], \\
    I_1 - I_2 &:= \{x-y\mid x\in I_1,\,y\in I_2\} = [a-d,b-c]
  \end{align*}
  for all $a,b,c,d\in\R$ with $a\le b$ and $c\le d$.
  
  \begin{lemma}\label{convolution_approximation}
    Let $g: \R \rightarrow \R$ be monotonic, bounded and positive,
    $k\in L^1_{\textnormal{loc}}(\R,\R)$ and $N>0$. If there are constants
    $\alpha_N,\beta_N>0$ such that
    \begin{itemize}
      \item $\abs{\int_I k(y)\dy} \le \alpha_N$ for all intervals $I\subseteq [-N,N]$,
      \item $\abs{\int_I k(y)\dy} \le \beta_N$ for all bounded intervals
            $I\subseteq \R\backslash[-N,N]$,
    \end{itemize}
    then
    \begin{equation*}
      \big(g*k\big)(x) \in \delta_N g(x-N) + \abs{g(x+N)-g(x-N)}[-\alpha_N, \alpha_N] + 4\norm{g}_\infty[-\beta_N, \beta_N],
    \end{equation*}
    where $\delta_N = \int_{[-N,N]} k(y)\dy$, holds for all $x\in\R$.
  \end{lemma}
  \proof Let $x\in\R$ and split the integral into three parts,
    \begin{equation*}
      \big(g\ast k\big)(x) = \int_{-N}^N g(x-y)k(y)\dy + \int_N^\infty g(x-y)k(y)\dy + \int_{-\infty}^{-N} g(x-y)k(y)\dy.
    \end{equation*}
    By the second mean value theorem of integral calculus \cite{hobson09} there exists a
    $B\in[-N,N]$ with
    \begin{align*}
      \int_{-N}^N g(x-y)k(y)\dy &= g(x+N)\int_{-N}^B k(y)\dy + g(x-N)\int_B^N k(y)\dy \\
      &= g(x+N)\int_{-N}^B k(y)\dy + g(x-N)\left(\delta_N - \int_{-N}^B k(y)\dy\right) \\
      &= \big(g(x+N)-g(x-N)\big)\int_{-N}^B k(y)\dy + \delta_N g(x-N) \\
      &\in \abs{g(x+N)-g(x-N)}[-\alpha_N,\alpha_N] + \delta_N g(x-N).
    \end{align*}
    Bounds for the other two integrals can also be found using the second mean value
    theorem: for every $M>N$ there exists a $Q_M\in[N,M]$ with
    \begin{align*}
      \int_N^\infty g(x-y)k(y)\dy &= \lim_{M\rightarrow\infty}\int_N^M g(x-y)k(y)\dy \\
      &= \lim_{M\rightarrow\infty} g(x-N)\int_N^{Q_M} k(y)\dy + g(x-M)\int_{Q_M}^M k(y)\dy \\
      & \in 2\norm{g}_\infty[-\beta_N, \beta_N]
    \end{align*}
    and similarly for the third integral.\qed
    
    \medskip
    If $k=\sinc$ is chosen as the convolution kernel in the previous lemma, the
    constants $\alpha_N$ and $\beta_N$ as well as a bound on $\delta_N$ can be found as
    follows.
    
    If $I = [c,d]\subset\R$ is a bounded interval, then
    \begin{equation*}
      \int_I k(y)\dy = \frac{1}{\pi}\int_{\pi c}^{\pi d} \frac{\sin(y)}{y}\dy = \frac{1}{\pi}\big(\Si(\pi d) - \Si(\pi c)\big),
    \end{equation*}
    where $\Si$ is the sine integral defined as $\Si(x) = \int_0^x {\sin(y)}/{y}\dy$
    for all $x\in\R$. It is well known that $\abs{\Si(x)}\le \Si(\pi) \le 2$
    for all $x\in\R$. Moreover, we have $\abs{\Si(x) - \frac{\pi}{2}} \le \frac{1}{x}$
    for all $x>0$, which can be shown with the Laplace transform. Because of the point
    symmetry of the sine integral it follows that
    $\abs{\Si(x) + \frac{\pi}{2}} \le \frac{1}{-x}$ for all $x<0$.
    
    Fix an arbitrary $N>0$. If $I = [c,d]\subseteq[-N,N]$, then
    \begin{equation*}
      \abs{\int_I k(y)\dy} = \frac{1}{\pi}\abs{\Si(\pi d) - \Si(\pi c)} \le \frac{1}{\pi}(2+2) = \frac{4}{\pi} =: \alpha_N,
    \end{equation*}
    On the other hand, if we consider an interval $I = [c,d]\subseteq(N,\infty)$, then
    \begin{align*}
      \abs{\int_I k(y)\dy} &= \frac{1}{\pi}\abs{\Si(\pi d) - \Si(\pi c)} = \frac{1}{\pi}\abs{\left(\Si(\pi d) - \frac{\pi}{2}\right) - \left(\Si(\pi c) - \frac{\pi}{2}\right)} \\
      &\le \frac{1}{\pi}\left(\frac{1}{\pi d} + \frac{1}{\pi c}\right) \le \frac{2}{\pi^2 N} =: \beta_N
    \end{align*}
    and similarly $\abs{\int_I f(y)\dy} \le \frac{2}{\pi^2 N}$ for
    $I = [c,d]\subseteq(-\infty,-N)$. The value of $\delta_N$ is bounded by
    \begin{equation*}
      \abs{\delta_N - 1} = \abs{\int_{[-N,N]}k(y)\dy - 1} = \abs{\frac{2}{\pi}\Si(\pi N) - 1} \le \frac{2}{\pi} \frac{1}{\pi N} = \frac{2}{\pi^2 N}.
    \end{equation*}
    
    \begin{lemma}\label{noncoercive_part1}
      The functional $f\mapsto\norm{f\star f}_{L^2}$ with $f\in L^2(I,\C)$  is not
      coercive for any interval $I\subseteq\R$.
    \end{lemma}
    \proof For $\delta>1$ let
    $g_\delta:\R \rightarrow \R$ with $g_\delta(x) = \frac{1}{\sqrt{x}}$ for all
    $x\in[1,\delta]$ and $g_\delta(x) = 0$ for all $x\in\R\backslash[1,\delta]$. The
    function $g_\delta$ is neither monotonic nor positive, but it can be written as the
    difference of two monotonic, bounded and positive functions $g_{\delta,1}$ and
    $g_{\delta,2}$, e.g.
    \begin{equation*}
      g_{\delta,1}(x) = \begin{cases}
                            1, & \text{if } x<1, \\
                            2, & \text{if } x\ge 1
                        \end{cases}
      \qquad\text{and}\qquad
      g_{\delta,2}(x) = \begin{cases}
                            1, & \text{if } x<1, \\
                            2-\frac{1}{\sqrt{x}}, & \text{if } 1\le x\le \delta, \\
                            2, & \text{if } x>\delta.
                        \end{cases}
    \end{equation*}
    Applying \cref{convolution_approximation} to $g_{\delta,j}$ for $j\in\{1,2\}$ and
    using $\alpha_n$ and $\beta_N$ as derived above as well as the bound on $\delta_N$
    yields bounds
    $\big(g_{\delta,j} * k\big)(x) \in [L_{\delta,j,N}(x), U_{\delta,j,N}(x)]$ for all
    $x\in\R$ and all $N>0$, where
    \begin{align*}
      L_{\delta,j,N}(x) &:= g_{\delta,j}(x-N)\left(1-\frac{2}{\pi^2N}\right) - \abs{g_{\delta,j}(x+N)-g_{\delta,j}(x-N)}\frac{4}{\pi} - \frac{16}{\pi^2N}, \\
      U_{\delta,j,N}(x) &:= g_{\delta,j}(x-N)\left(1+\frac{2}{\pi^2N}\right) + \abs{g_{\delta,j}(x+N)-g_{\delta,j}(x-N)}\frac{4}{\pi} + \frac{16}{\pi^2N}.
    \end{align*}
    Therefore, $g_{\delta}*k$ is bounded by
    \begin{align*}
      \big(g_\delta * k\big)(x) &= \big(g_{\delta,1} * k\big)(x) - \big(g_{\delta,2} * k\big)(x) \\
      &\in [L_{\delta,1,N}(x) - U_{\delta,2,N}(x), U_{\delta,1,N}(x) - L_{\delta,2,N}(x)]
    \end{align*}
    for all $x\in\R$ and all $N>0$. In particular, choosing $x=2N$ yields the lower
    bounds
    \begin{align*}
      \big(g_\delta * k\big)(2N) \ge
      \begin{cases}
        \left(1 - \frac{12 - 4\sqrt{3}}{3\pi}\right)\frac{1}{\sqrt{N}} - \frac{40}{\pi^2}\frac{1}{N} + \frac{2}{\pi^2}\frac{1}{N\sqrt{N}}, & \text{if } N \in \left[1,\frac{\delta}{3}\right], \\
        \left(1 - \frac{4}{\pi}\right)\frac{1}{\sqrt{N}} - \frac{40}{\pi^2}\frac{1}{N} + \frac{2}{\pi^2}\frac{1}{N\sqrt{N}}, & \text{if } N \in \left(\frac{\delta}{3},\delta\right], \\
        -\frac{40}{\pi^2N}, & \text{if } N \in \left(\delta,\infty\right)
      \end{cases}
    \end{align*}
    and the upper bounds
    \begin{align*}
      \big(g_\delta * k\big)(2N) \le
      \begin{cases}
        \left(1 + \frac{12 - 4\sqrt{3}}{3\pi}\right)\frac{1}{\sqrt{N}} + \frac{40}{\pi^2}\frac{1}{N} - \frac{2}{\pi^2}\frac{1}{N\sqrt{N}}, & \text{if } N \in \left[1,\frac{\delta}{3}\right], \\
        \left(1 + \frac{4}{\pi}\right)\frac{1}{\sqrt{N}} + \frac{40}{\pi^2}\frac{1}{N} - \frac{2}{\pi^2}\frac{1}{N\sqrt{N}}, & \text{if } N \in \left(\frac{\delta}{3},\delta\right], \\
        \frac{40}{\pi^2N}, & \text{if } N \in \left(\delta,\infty\right)
      \end{cases}
    \end{align*}
    for all $N>1$. A bound on $\big(g_\delta * k\big)(x)$ for negative values of $x$ can
    be found similarly by choosing $x=-N$, which yields
    $\big(g_\delta * k\big)(-N) \in \left[\frac{-36}{\pi^2 N}, \frac{36}{\pi^2 N}\right]$
    for all $N>0$. Furthermore, $g_\delta * k$ is uniformly bounded by
    \begin{equation*}
      \norm{g_\delta * k}_\infty \le \norm{g_\delta}_{L^4}\norm{k}_{L^{4/3}} \le \norm{g_\infty}_{L^4}\norm{k}_{L^{4/3}} < \infty
    \end{equation*}
    for all $\delta>1$ by Hölder's inequality. Using these bounds and the fact that the
    coefficient $1 - \frac{12 - 4\sqrt{3}}{3\pi}$ of $\frac{1}{\sqrt{N}}$ in the lower
    bound of $\big(g_\delta*k\big)(2N)$ is positive, it follows that 
    \begin{equation*}
      \lim_{\delta\rightarrow\infty}\norm{g_\delta * k}_{L^2} = \infty \qquad\text{and}\qquad
      \lim_{\delta\rightarrow\infty}\norm{(g_\delta * k)^2}_{L^2} < \infty.
    \end{equation*}
    
    Without loss of generality we assume that $I=\left[-\frac{1}{2},\frac{1}{2}\right]$.
    If $f_\delta = (\mathcal{F}g_\delta)\chi_{[-\frac{1}{2},\frac{1}{2}]}$ for all
    $\delta>1$, then $f_\delta\in L^2\left(\left[-\frac{1}{2}, \frac{1}{2}\right], \C\right)$ and
    \begin{align*}
      \lim_{\delta\rightarrow\infty}\norm{f_\delta}_{L^2} &= \lim_{\delta\rightarrow\infty}\norm{\mathcal{F}^{-1}f_\delta}_{L^2} = \lim_{\delta\rightarrow\infty}\norm{g_\delta \ast k}_{L^2} = \infty, \\
      \lim_{\delta\rightarrow\infty}\norm{f_\delta \star f_\delta}_{L^2} &= \lim_{\delta\rightarrow\infty}\norm{\big(\mathcal{F}^{-1}f_\delta\big)^2}_{L^2} = \lim_{\delta\rightarrow\infty}\norm{(g_\delta \ast k)^2}_{L^2} < \infty. \tag*{$\qed$}
    \end{align*}
    
    \begin{corollary}\label{noncoercive_part2}
      The functional $f\mapsto\norm{f\star f}_{L^2}$ with $f\in L^2(U,\C)$ is not
      coercive for any $U\subseteq\R^d$ with non-empty interior.
    \end{corollary}
    \proof Since the interior of $U$ is non-empty, there are $a,b\in\R$ with $a<b$ such
      that $[a,b]^d\subseteq U$. By \cref{noncoercive_part1} there is a sequence
      $(f_n)_{n\in\N}$ of functions $f_n\in L^2([a,b],\C)$ with
      $\lim_{n\rightarrow\infty} \norm{f_n}_{L^2} = \infty$ and
      $\lim_{n\rightarrow\infty} \norm{f_n\star f_n}_{L^2} = 1$. Consider the sequence
      $(h_n)_{n\in\N}$ given by
      \begin{equation*}
        h_n: \R^d \rightarrow \C, \quad x \mapsto \begin{cases} f_n(x_1), & \text{if } x \in [a,b]^d, \\ 0, & \text{otherwise.} \end{cases}
      \end{equation*}
      Then
      \begin{equation*}
        \norm{h_n}_{L^2}^2 = \int_{[a,b]^d} \abs{f_n(x_1)}^2 \dx = (b-a)^{d-1} \norm{f_n}_{L^2}^2 \;\underset{n\rightarrow\infty}{\longrightarrow}\; \infty
      \end{equation*}
      and
      \begin{align*}
        \abs{\big(h_n\star h_n\big)(x)} &= \abs{\int_{\R^d} h_n(y)h_n(x+y)\dy} = \abs{\int_{[a,b]^d \cap ([a,b]^d - x)} f_n(y_1)f_n(y_1+x_1)\dy} \\
        &\le \begin{cases}
               (b-a)^{d-1}\abs{\big(f_n\star f_n\big)(x_1)}, & \text{if } x \in B_{b-a}^\infty(0), \\
               0, & \text{otherwise}
             \end{cases}
      \end{align*}
      for all $x\in\R^d$, where $B_{b-a}^\infty(0) = \left\{x\in\R^d: \norm{x}_\infty < b-a\right\}$. This implies
      \begin{align*}
        \norm{h_n\star h_n}_{L^2}^2 &\le (b-a)^{2d-2}\int_{B_{b-a}^\infty(0)}\abs{\big(f_n\star f_n\big)(x_1)}^2\dx \\
        &\le (b-a)^{2d-2} (b-a)^{d-1} \int_{[-(b-a),b-a]} \abs{\big(f_n\star f_n\big)(x_1)}^2\dx \\
        &= (b-a)^{3d-3} \norm{f_n\star f_n}_{L^2}^2 \;\underset{n\rightarrow\infty}{\longrightarrow}\; (b-a)^{3d-3} < \infty. \tag*{\qed}
      \end{align*}
    
    An immediate consequence of \cref{noncoercive_part2} is that the functional
    \begin{equation*}
      L^2(A,\C) \mapsto \R, \qquad \Psi \mapsto \norm{(\Psi p_Z a) \star (\Psi p_Z a) - G}_{L^2}^2
    \end{equation*}
    is not coercive for any $G\in L^2(\R^2, \C)$, since $p_Z(v)\neq 0$ for all $v\in\R^2$
    and $A$ is the set-theoretic support of $a$. This functional is equal to $E$ in the
    case of perfectly coherent illumination and only one experimental input image.
  
    \section*{References}
    \bibliographystyle{plain}
    \bibliography{summary}
  
\end{document}